\documentclass{article}

\def\bewname{Beweis}
\def\definname{Definition}
\def\thmname{Theorem}
\def\satzname{Satz}
\def\lemmaname{Lemma}
\def\folgname{Folgerung}
\def\bemname{Bemerkung}
\def\bspname{Beispiel}

\newcounter{beh}[section]

\def\thebeh{\thesection.\arabic{beh}}

\newenvironment{bew}[1]{{\sc \bewname} #1:}{\ \hspace*{\fill}
\rule{1ex}{1ex}}

\newenvironment{defin}{\refstepcounter{beh}{\sc \definname\
\thebeh.} }{}

\newenvironment{thm}{\refstepcounter{beh}{\sc \thmname\
\thebeh.} \em }{\rm}

\newenvironment{satz}{\refstepcounter{beh}{\sc \satzname\
\thebeh.} \em }{\rm}

\newenvironment{lemma}{\refstepcounter{beh}{\sc \lemmaname\
\thebeh.} \em }{}

\newenvironment{folg}{\refstepcounter{beh}{\sc \folgname\
\thebeh.} \em }{}

\newcommand{\Rda}{\hat{R}}

\newcommand{\Rdam}{{\hat{R}^{-1}}}

\newcommand{\Gammalinv}{{_{\mathrm{inv}}\Gamma}}

\newcommand{\linv}[1]{{_{\mathrm{inv}}{#1}}}

\newcommand{\ad}{\mathrm{ad\,}}
\newcommand{\s}{\mathfrak{s}}
\newcommand{\cont}{^{\mathrm{c}}}

\newcommand{\dif}{\mathrm{d}}

\def\bewname{Proof}
\def\definname{Definition}
\def\thmname{Theorem}
\def\satzname{Proposition}
\def\lemmaname{Lemma}
\def\folgname{Corollary}
\def\bemname{Remark}
\def\bspname{Example}

\usepackage{amstex}
\newcommand{\z}{p}
\newcommand{\adR}{\ad\!_R}

\def\gwP{P\kern-1.5ex\lower-1.5ex\hbox{\smash{\scriptsize $\rightharpoonup$}}}
\def\gwomega{\omega\kern-1.6ex\lower-.92ex\hbox{\smash{\scriptsize $\rightharpoonup$}}}
\def\gwY{Y\kern-1.8ex\lower-1.5ex\hbox{\smash{\scriptsize $\rightharpoonup$}}}

\begin{document}

\begin{center}
\LARGE Levi-Civita Connections on the Quantum Groups $SL_{q}(N)$,\\
$O_{q}(N)$ and $Sp_{q}(N)$
\end{center}

\vspace{1.8cm}

{\bf Istv{\'a}n~Heckenberger, Konrad~Schm{\"u}dgen,}
Universit{\"a}t Leipzig, Fakult{\"a}t f{\"u}r Mathematik und
Informatik und NTZ, Augustusplatz 10, 04109 Leipzig, Germany.\\
{\scriptsize e-mail:\,istvan@@aix550.informatik.uni-leipzig.de,
schmuedgen@@mathematik.uni-leipzig.d400.de}

\setcounter{section}{-1}

\begin{abstract}
For bicovariant differential calculi on quantum groups various notions
on connections and metrics (bicovariant connections, invariant
metrics, the compatibility of a connection with a metric, Levi-Civita
connections) are introduced and studied. It is proved that for the
bicovariant differential calculi on $SL_{q}(N)$, $O_{q}(N)$ and
$Sp_{q}(N)$ from the classification in \cite{a-SchSch1} there
exist unique Levi-Civita connections.
\end{abstract}

\section{Introduction} \label{sec-einfuehrung}

The seminal work of S.\,L.\ Woronowicz \cite{a-Woro2} was the starting point
to study non-commutative bicovariant differential calculi on quantum
groups (Hopf algebras). Woronowicz has developed a general theory of
such calculi which in many aspects can be considered as a non-commutative
version of the classical Lie group theory. Bicovariant differential
calculi on the quantum matrix groups $SL_{q}(N)$, $O_{q}(N)$ and
$Sp_{q}(N)$ have been classified (under natural assumptions) in two recent
papers \cite{a-SchSch1} and \cite{a-SchSch2}.
An outcome of this classification is that except for
finitely many values of $q$ there are precisely $2N$ such calculi on
$SL_{q}(N)$ for $N\geq 3$ and two on $O_{q}(N)$ and $Sp_{q}(N)$ for
$N\geq 4$. See Section \ref{sec-jurcokalkule} for a brief review.
It is clear that these calculi are basic tools of non-commutative
geometry on the corresponding quantum groups.

The aim of this paper is to define and to study invariant metrics and
Levi-Civita connections for the bicovariant differential calculi on
$SL_{q}(N)$, $O_{q}(N)$ and $Sp_{q}(N)$. The Killing metric of a
compact Lie group and the Levi-Civita connection  of a Riemannian
manifold are fundamental notions of (commutative) differential geometry,
so it seems that extensions of these concepts to the bicovariant
differential calculi are necessary steps toward the development of
non-commutative differential geometry on quantum groups. It turns out
that these generalizations are by no means straightforward and that
phenomena occur which are absent in classical differential geometry.
We briefly mention some of these new features. Firstly, no ad-invariant
metric for these bicovariant differential calculi is symmetric in the
usual sense if $q$ is not a root of unity. However, the ad-invariant
metrics are symmetric with respect to the corresponding braiding, see
Corollary \ref{f-symmmetrik} below. In our opinion, this fits nicely
into the concepts of the braided geometry, see \cite{b-Majid1} and the
references therein. Secondly, ad-invariant metrics for
$SL_{q}(N)$ resp.\ $O_{q}(N)$, $Sp_{q}(N)$ depend on two resp.\ three
complex parameters. Roughly speaking, this stems from the fact that
the dimension of these bicovariant calculi is $N^{2}$ rather than the
dimension of the corresponding Lie groups as in case of the classical
differential calculus. Thirdly, if we generalize the notion of a
Levi-Civita connection in obvious straightforward manner, then there
exist many Levi-Civita connections for a given metric. For our
bicovariant differential calculi they depend on a number of free
complex parameters, see e.\,g.\ appendix \ref{sec-LCZ2}.

The main purpose of this paper is to propose notions of compatibility of
a connection and an ad-invariant metric (see formula
(\ref{eq-adinvmetrik}) and Definition \ref{d-vertrmet} below)
and of a Levi-Civita connection
(see Definition \ref{d-LCZ}) which ensure that the important result
in classical differential geometry of uniqueness of the Levi-Civita
connection remains valid in the present setting. If the braiding map
$\sigma$ is the flip operator of the tensor product as in case of the
``ordinary'' differential calculus on compact Lie groups, our
definitions give just the corresponding classical concepts.

This paper is organized as follows. In Section \ref{sec-jurcokalkule}
we briefly recall the bicovariant differential calculi on $SL_{q}(N)$,
$O_{q}(N)$ and $Sp_{q}(N)$ studied in
\cite{a-SchSch1} and we collect some basic facts needed later. In Section
\ref{sec-metrik} we define metrics and invariant metrics. The invariant
metrics for these calculi and their restrictions to the invariant
subspaces for the right coaction and the right adjoint action,
respectively, are determined. For some of the calculi the quantum Lie
algebra contains the corresponding classical Lie algebra as
ad-invariant subspace when $q\to1$. In these cases the limits of the
invariant metrics exist and their restrictions to these subspaces give
multiples of the Killing forms. In Appendix \ref{sec-Rossoform} we
show how a variant of the Rosso form of the quantized enveloping
algebra ${\cal U}_{q}(sl(N))$ gives an ad-invariant metric for
$SL_{q}(N)$. Section \ref{sec-zusammenhang} is
concerned with connections. After reviewing some generally known
definitions on connections (see e.\,g.\ \cite{b-Connes1}) we define
bicovariant connections, the compatibility of a connection with a metric
and Levi-Civita connections and we discuss some of their properties. The
main results of this paper (Theorems \ref{t-lcza} and \ref{t-lczbcd}) are
stated and proved in Sections \ref{sec-lcza} and \ref{sec-lczbcd},
respectively. They assert that for each of the bicovariant differential
calculi on $SL_{q}(N)$, $O_{q}(N)$ and $Sp_{q}(N)$ described in Section
\ref{sec-jurcokalkule} there is precisely one Levi-Civita connection.
Moreover, it is shown that these Levi-Civita connections admit limits
when $q\to1$ in an appropriate way. In
appendix \ref{sec-LCZ2} of this paper we show that if we define
Levi-Civita connections by taking the ``usual'' compatibility with a
metric, then the set of Levi-Civita connections for $SL_{q}(N)$ depends
on three free parameters.

Let us fix some notation and assumptions which are needed in the sequel.
Throughout we use Sweedler's notation $\Delta(a)=a_{(1)}\otimes a_{(2)}$
and the Einstein convention to sum over repeated indices. The antipode
of a Hopf algebra is denoted by $\kappa$ and the counit by
$\varepsilon$. Let $\mbox{Mor}(v,w)$ the space of intertwiners
of corepresentations $v$ and $w$ and let $\mbox{Mor}(v):=
\mbox{Mor}(v,v)$. We use the definitions of the quantum groups
$SL_{q}(N)$, $O_{q}(N)$ and $Sp_{q}(N)$ and their basic properties
established in \cite{a-FadResTak1}. Let $u=(u^{i}_{j})$ denote the
corresponding fundamental representation and let $u\cont$ be the
contragredient representation of $u$. In case where ${\cal A}$ is
$O_{q}(N)$ or $Sp_{q}(N)$ let $C=(C{^{i}_{j}})$ denote the matrix of
the metric of ${\cal A}$ and $B$ the inverse matrix of $C$.
We abbreviate $Q=q-q^{-1}$.

Also we shall freely use the general theory  of bicovariant
differential calculi developed in \cite{a-Woro2}. The abbreviation FODC
means a first order differential calculus. As noted above, our main
intention is to study the bicovariant differential calculi
$\Gamma_{\pm,z}$ (see Section \ref{sec-jurcokalkule}) which occured
in the classification of \cite{a-SchSch1,a-SchSch2}. But the
corresponding concepts and general facts apply to an arbitrary
bicovariant differential calculus over a Hopf algebra. In order to
avoid confusion, let us adopt the following notations:
$\tilde{\Gamma}$ is always a general bicovariant FODC over a Hopf
algebra ${\cal A}$, ${\cal A}^{0}$ is the Hopf dual of ${\cal A}$,
$\tilde{\cal X}$ is the quantum Lie algebra of $\tilde{\Gamma}$,
$\{\eta_{i}\,|\,i\in{\cal I}\}$ is a finite vector space basis of
$\linv{\tilde{\Gamma}}$, $\{\chi_{i}\,|\,i\in{\cal I}\}$ is the
corresponding dual basis of $\tilde{\cal X}$ and $\tilde{\sigma}$ is
the braiding of $\tilde{\Gamma}\otimes_{\cal A}\tilde{\Gamma}$ as
defined by Proposition 3.1 in \cite{a-Woro2}.

In this paper we suppose that the deformation parameter $q$ {\bf is not a
root of unity} and $q\not= 0$. Then, roughly speaking, the representation
theory of ${\cal A}$ is similar to the classical case
\cite{a-Lusztig1,a-Rosso2}. We shall need this assumption only in order to
ensure that the decompositions of certain tensor product representations
of $u$ and $u\cont$ can be labelled by Young tableaus similar to the
classical case. The corresponding results are the Lemmas
\ref{l-invmetrika}, \ref{l-invmetrikbcd}, \ref{l-zusa} and
\ref{l-zusbcd}. All other considerations are valid without this
assumption.

\section{Review of some facts on bicovariant differential calculi
on quantum groups $SL_{q}(N)$, $O_{q}(N)$ and $Sp_{q}(N)$}
\label{sec-jurcokalkule}

Firstly we repeat the construction of bicovariant differential calculi,
see \cite{a-SchSch1} for some missing details of proofs in the following
discussion.

Let $z$ be a nonzero complex number. We assume that $z^{N}=q$ for
${\cal A}=SL_{q}(N)$ and that $z^{2}=1$ for ${\cal A}=O_{q}(N)$ and
${\cal A}=Sp_{q}(N)$.
Let $L{^{\pm}_{z}}=(l{^{\pm}_{z}}{^{i}_{j}})$ denote the $N\times N$
matrix of linear functionals $l{^{\pm}_{z}}{^{i}_{j}}$ on ${\cal A}$
defined in \cite{a-FadResTak1}, Section 2, by taking the matrix
$z^{-1}P\Rda$ as $R$. By definition, we then have
\begin{equation} \label{eq-lfunktionale}
l{^{+}_{z}}{^{i}_{j}}(u{^{n}_{m}})=z^{-1}\Rda{^{in}_{mj}}\quad
\mbox{ and }\quad l{^{-}_{z}}{^{i}_{j}}(u{^{n}_{m}})=z\Rdam{^{in}_{mj}}.
\end{equation}
Let $D_{i}:=q^{2i}$ for ${\cal A}=SL_{q}(N)$ and
$D_{i}:=(C(C^{-1})^{\rm t}){^{i}_{i}}$ for ${\cal A}=O_{q}(N)$,
$Sp_{q}(N)$, where $C$ is the matrix given by the metric of $O_{q}(N)$
and $Sp_{q}(N)$, cf.\ \cite{a-FadResTak1}, Sect.\ 1. Then we have
$\kappa^{2}(u^{i}_{j})=D_{i}u^{i}_{j}D_{j}^{-1}$ (no summation).

There are $2N$ bicovariant FODC $\Gamma_{\pm,z}$, $z^N=q^2$ on
$SL_q(N)$, $N\geq3$, and 2 bicovariant FODC $\Gamma_+=\Gamma_{+,1}$
and $\Gamma_-=\Gamma_{+,-1}$ on $O_q(N)$, $Sp_q(N)$ and $SL_q(2)$.
Except for the quantum group $O_q(3)$, these FODC exhaust the bicovariant
first order differential calculi which appeared in the classification of
\cite{a-SchSch1,a-SchSch2}. They are the objects of our study in the paper.

In what follows we assume that $\Gamma$ always denotes such a FODC
$\Gamma_{\pm,z}$. The FODC $\Gamma_{\pm,z}$ on ${\cal A}$ is given by
\[\dif a=\sum_{i,j=1}^{N}(\chi_{ij}*a)\eta_{ij},\quad a\in{\cal A},\]
where
\begin{equation} \label{eq-chifunkt}
\chi_{ij}=\sum_{n=1}^ND_n^{-1}l{^\pm_{z'}}{^n_i}\kappa(l{^\mp_{z''}}{^j_n})
-D_i^{-1}\delta_{ij}\varepsilon,\quad i,j=1,\ldots,N
\end{equation}
and $\{\eta_{ij}\,|\,i,j=1,\ldots,N\}$ is a basis of the space
$\linv{(\Gamma_{\pm,z})}$ of left-invariant elements of
$\Gamma_{\pm,z}$. The
right and left ${\cal A}$-module operations satisfy the equations
$\eta_{ij}a=(f{^{ij}_{mn}}*a)\eta_{mn}$, $a\in{\cal A}$, and the right
coaction of ${\cal A}$ on this basis is given by $\Delta_{R}(\eta_{ij}
)=\eta_{mn}\otimes v{^{mn}_{ij}}$, $i,j=1,\ldots,N$, where
\begin{equation} \label{eq-bimodulvertauschung}
f{^{ij}_{mn}}=l{^{\pm}_{z'}}{^{i}_{m}}\kappa(l{^{\mp}_{z''}}{^{n}_{j}}),
\end{equation}
\begin{equation} \label{eq-rechtskowirkung}
v{^{mn}_{ij}}=(u\cont u)^{mn}_{ij}=\kappa(u^{i}_{m})u^{n}_{j}.
\end{equation}
The linear span ${\cal X}$ of the linear functionals $\chi_{ij}$,
$i,j=1,\ldots,N$, equipped with the bracket $[\cdot,\cdot]:{\cal X}
\times{\cal X}\rightarrow{\cal X}$ defined in \cite{a-Woro2}, Section 5,
is called the {\bf quantum Lie algebra} of the bicovariant FODC
$\Gamma_{\pm,z}$.

There is a duality $<\cdot,\cdot>$ between $\linv{\tilde{\Gamma}}$ and
$\tilde{\cal X}$ given by
\begin{equation} \label{eq-dualitaet}
\left<\sum_{i}a_{i}\eta_{i},\sum_{j}\chi_{j}b_{j}\right> =
\sum_{k}a_{k}b_{k}
\end{equation}
for all $a_{k},b_{k}\in\mathbb{C}$. This definition extends to a map
$<\cdot,\cdot>: A\otimes\linv{\tilde{\Gamma}}\times\tilde{\cal X}
\otimes B\rightarrow A\otimes_{\cal A}B$ for linear subspaces $A$ of
$\tilde{\Gamma}^{\wedge}$ resp.\ $\tilde{\cal X}\otimes{\cal A}$ and
linear subspaces $B$ of $\tilde{\Gamma}^{\wedge}$ by taking
$a_{k}\in A$ and $b_{k}\in B$ in (\ref{eq-dualitaet}) for $k\in{\cal I}$.

Let $\sigma$ denote the braiding map of $\Gamma_{\pm,z}\otimes_{\cal A}
\Gamma_{\pm,z}$ (see Proposition 3.1 in \cite{a-Woro2}) and let
$\sigma^{mnrs}_{ijkl}$ be the matrix coefficients of $\sigma$
with respect to the basis $\{\eta_{ij}\otimes_{\cal A}\eta_{kl}\}$
of $\linv{(\Gamma_{\pm,z}\otimes_{\cal A}\Gamma_{\pm,z})}$,
i.\,e.\ $\sigma(\eta_{ij}\otimes\eta_{kl})=\sigma{^{mnrs}_{ijkl}}
\eta_{mn}\otimes\eta_{rs}.$

\begin{lemma} \label{l-sigmamatrix}
$\sigma{^{mnrs}_{ijkl}}=D_{k}D_{x}^{-1}\Rda^{\mp}{^{pk}_{xj}}
\Rda^{\pm}{^{tn}_{yr}}\Rda^{\mp}{^{xi}_{tm}}\Rda^{\pm}{^{ys}_{pl}},$
where the upper and lower signs refer to $\Gamma_{+,z}$ and
$\Gamma_{-,z}$, respectively.
\end{lemma}

\begin{bew}{}
We carry out the proof for $\Gamma_{+,z}$. {}From
(\ref{eq-bimodulvertauschung}) we see that $f{^{tp}_{rs}}=
l{^{+}_{z'}}{^{t}_{r}}\kappa(l{^{-}_{z''}}{^{s}_{p}})$ are the linear
functionals of \cite{a-Woro2}, Theorem 2.1, for the FODC $\Gamma_{+,z}$.
(\ref{eq-lfunktionale}) yields $f{^{tp}_{rs}}(u^{n}_{l})=
z^{-1}\Rda{^{tn}_{yr}}\Rda{^{ys}_{pl}}$ and so $\delta{^{i}_{r}}
\delta{^{j}_{s}}\delta{^{n}_{m}}=f{^{ij}_{rs}}(\kappa(\kappa(u^{n}_{k})
u^{k}_{m}))=f{^{ij}_{tp}}(\kappa(u^{k}_{m}))
f{^{tp}_{rs}}(\kappa^{2}(u^{n}_{k}))$.
{}From $\kappa^{2}(u^{n}_{k})=D_{n}D_{k}^{-1}u^{n}_{k}$ we easily derive
$\delta{^{i}_{r}}\delta{^{j}_{s}}\delta{^{n}_{m}}=
f{^{ij}_{tp}}(\kappa(u^{k}_{m}))z^{-1}D_{n}D_{k}^{-1}\Rda{^{tn}_{xr}}
\Rda{^{xs}_{pk}}$, so the latter yields
$f{^{ij}_{tp}}(\kappa(u^{k}_{m}))=zD_{k}D_{x}^{-1}\Rdam{^{pk}_{xj}}
\Rdam{^{xi}_{tm}}$.
By the general theory (cf.\ formula (3.15) in \cite{a-Woro2}) and
(\ref{eq-rechtskowirkung}), we have $\sigma{^{mnrs}_{ijkl}}=
f{^{ij}_{rs}}((u\cont)^{m}_{k}u^{n}_{l})=
f{^{ij}_{tp}}(\kappa(u^{k}_{m}))f{^{tp}_{rs}}(u^{n}_{l})$ from which the
above formula follows.
\end{bew}

Recall that the matrix $\Rda$ satisfies the Hecke relation
$(\Rda-qI)(\Rda+q^{-1}I)=0$ for ${\cal A}=SL_{q}(N)$ and the cubic
equation $(\Rda-qI)(\Rda+q^{-1}I)(\Rda-\epsilon q^{\epsilon-N}I)=0$
for ${\cal A}=O_{q}(N)$, $Sp_{q}(N)$, where $\epsilon=1$ for
${\cal A}=O_{q}(N)$ and $\epsilon=-1$ for ${\cal A}=Sp_{q}(N)$. {}From
these equations and Lemma \ref{l-sigmamatrix} it follows that
\[(\sigma-I)(\sigma+q^{-2}I)(\sigma+q^{2}I)=0\]
for ${\cal A}=SL_{q}(N)$ and
\[(\sigma-I)(\sigma+q^{-2}I)(\sigma+q^{2}I)
(\sigma-\epsilon q^{N-\epsilon+1}I)(\sigma-\epsilon q^{\epsilon-N-1}I)
(\sigma+\epsilon q^{N-\epsilon-1}I)(\sigma+\epsilon q^{\epsilon-N+1}I)=0\]
for ${\cal A}=O_{q}(N)$, $Sp_{q}(N)$ (cf.\ \cite{a-CSchWW1}).

The above formulas show that the set of for $q>0$ positive
eigenvalues of $\sigma$ is $\{1\}$ for $SL_{q}(N)$ and
$\{1,q^{N},q^{-N}\}$ for $O_{q}(N)$ and $Sp_{q}(N)$.

In the following we use the abbreviations $\z:=\epsilon q^{N-\epsilon}$,
$\s:=1+Q^{-1}(\z-\z^{-1})$ for $O_{q}(N)$ and $Sp_{q}(N)$ and
$\s:=\sum_{i=1}^{N}q^{-2i}$ for $SL_{q}(N)$.

Let $\omega_{ij}:=\kappa(u^{i}_{n})\dif u^{n}_{j}$, $i,j=1,\ldots,N$.
In \cite{a-SchSch1} and \cite{a-SchSch2} it is proved that if $q$ is not a
root of unity and apart from finitely many other values of $q$ the set
$\{\omega_{ij}\,|\,i,j=1,\ldots,N\}$ is a basis of the vector space
$\Gammalinv$. Let $\{X_{ij}\,|\,i,j=1,\ldots,N\}$ be the corresponding
dual basis of ${\cal X}$ with respect to the duality
(\ref{eq-dualitaet}). We call the sets $\{\omega_{ij}\}$ and
$\{X_{ij}\}$ {\bf standard bases} of $\Gammalinv$ and ${\cal X}$,
respectively. For this basis of $\Gammalinv$ the right coaction also
fulfills $\Delta_{R}(\omega_{ij})=\omega_{kl}\otimes v^{kl}_{ij}$
(cf.\ (\ref{eq-rechtskowirkung})).

Let us briefly return to a general bicovariant FODC $\tilde{\Gamma}$
over a Hopf algebra ${\cal A}$ such that
$\mathrm{dim}\linv{\tilde{\Gamma}}<\infty$. Then the quantum Lie algebra
$\tilde{\cal X}$ of $\tilde{\Gamma}$ is contained in the Hopf dual
${\cal A}^0$ and the bracket $[\cdot,\cdot]$ of $\tilde{\cal X}$ can be
written as $[x,y]=\adR y(x)$. Here $\adR$ denotes the right adjoint
action of ${\cal A}^0$ given by $\adR f(g)=\kappa(f_{(1)})gf_{(2)}$
for $f,g\in{\cal A}^0$. Moreover, we have $[x,f]:=\adR f(x)\in
\tilde{\cal X}$ for all $f\in{\cal A}^0$ and $x\in\tilde{\cal X}$.
A linear subspace ${\cal Y}$ of $\tilde{\cal X}$ is called ad-invariant
if $[y,f]\in{\cal Y}$ for all $y\in{\cal Y}$ and $f\in{\cal A}^0$.

Next we are looking for ad-invariant subspaces of $\tilde{\cal X}$.
For this we need the following simple

\begin{lemma} \label{l-adinvsubalg}
Suppose that $A=(A^{i}_{j})\in\mbox{Mor}(\tilde{v})$, i.\,e.\
$A^{i}_{j}\tilde{v}^{j}_{k}=\tilde{v}^{i}_{l}A^{l}_{k}$ for all
$i,k\in{\cal I}$. Then $\mbox{\rm im}\,A=\mbox{\rm lin}\{A^{i}_{j}
\eta_{i}\,|\,j\in{\cal I}\}$ is a $\Delta_{R}$-invariant subspace of
$\linv{\tilde{\Gamma}}$ and $\mbox{\rm im}\,A^{\rm t}=
\mbox{\rm lin}\{A^{i}_{j}\chi_{j}\,|\,i\in{\cal I}\}$ is
an ad-invariant subspace of $\tilde{\cal X}$.
\end{lemma}

\begin{bew}{}
Since $A\in\mbox{Mor}(\tilde{v})$, we have $\Delta_{R}(A^{i}_{j}
\eta_{i})=A^{i}_{j}\eta_{k}\otimes\tilde{v}^{k}_{i}=A^{k}_{l}\eta_{k}
\otimes\tilde{v}^{l}_{j}$, so that $\Delta_{R}(\mbox{im}\,A)\subset
\mbox{im}\,A\otimes{\cal A}$.
From the general theory \cite{a-Woro2} we easily derive that
$[\chi_{i},f]=f(\tilde{v}^{i}_{k})\chi_{k}$ for all $f\in{\cal A}^0$.
Therefore, we compute $\adR f(A^{i}_{k}\chi_{k})=A^{i}_{k}
f(\tilde{v}^{k}_{l})\chi_{l}=f(\tilde{v}^{i}_{j})A^{j}_{l}\chi_{l}$.
That is, $\adR f(\mbox{im}\,A^{\rm t})\subset\mbox{im}\,A^{\rm t}$ for
all $f\in{\cal A}^{0}$.
\end{bew}

For $SL_{q}(N)$ the projections $P_{0}$, $P_{1}:\Gammalinv\rightarrow
\Gammalinv$ defined by $P_{0}{^{ij}_{kl}}=\frac{1}{\s}q^{-2i}
\delta^{ij}\delta_{kl}$, $P_{1}{^{ij}_{kl}}=\delta^{i}_{k}
\delta^{j}_{l}-P_{0}{^{ij}_{kl}}$ belong to $\mbox{Mor}(u\cont\otimes
u)$. Hence the subspaces $\Upsilon^{0}=\mbox{im}\,P_{0}=\mbox{lin}\{
\omega^{0}:=P_{0}(\omega_{11})=P_{0}{^{kl}_{11}}\omega_{kl}\}$ and
$\Upsilon^{1}=\mbox{im}\,P_{1}=\mbox{lin}\{\omega^{1}_{ij}:=
P_{1}(\omega_{ij})=P_{1}{^{kl}_{ij}}\omega_{kl}\}$ of $\Gammalinv$ are
$\Delta_{R}$-invariant. The corresponding ad-invariant subspaces of
${\cal X}$ are ${\cal Y}^{0}=\mbox{im}\,P_{0}^{\rm t}=\mbox{lin}\{
Y^{0}:=\frac{1}{\s}\sum_{k}X_{kk}\}$ and ${\cal Y}^{1}=\mbox{im}\,
P_{1}^{\rm t}=\mbox{lin}\{Y^{1}_{ij}:=P_{1}{^{ij}_{kl}}X_{kl}=X_{ij}-
q^{-2i}\delta^{ij}Y^{0}\}$ respectively. Since $\{\omega_{ij}\}$ and
$\{X_{ij}\}$ are dual bases, we also have $\dif a=\sum_{ij}(X_{ij}*a)
\omega_{ij}$ for all $a\in{\cal A}$. Because of $P_{0}+P_{1}=
\mbox{id}$ and $P_{0}P_{1}=P_{1}P_{0}=0$, the latter leads to the
formula
\[\dif a=\sum_{i,j=1}^{N}(Y^{1}_{ij}*a)\omega^{1}_{ij}+\s(Y^{0}*a)
\omega^{0},\quad a\in{\cal A}.\]

Now we turn to the quantum groups $O_{q}(N)$ and $Sp_{q}(N)$. Then the
projections $P_{0}$, $P_{+}$, $P_{-}:\Gammalinv\rightarrow\Gammalinv$
given by
\[P_{0}{^{ij}_{kl}}=\frac{1}{\s}B^{\rm t}{^{i}_{m}}C^{mj}\delta_{kl},
\quad
P_{+}{^{ij}_{kl}}=\frac{1}{q+q^{-1}}(q^{-1}\delta^{i}_{k}\delta^{j}_{l}
+B^{\rm t}{^{i}_{m}}\Rda{^{mj}_{nl}}C^{\rm t}{^{n}_{k}}-
(q^{-1}+\z^{-1})P_{0}{^{ij}_{kl}}),\]
\[P_{-}{^{ij}_{kl}}=\frac{1}{q+q^{-1}}(q\delta^{i}_{k}\delta^{j}_{l}-
B^{\rm t}{^{i}_{m}}\Rda{^{mj}_{nl}}C^{\rm t}{^{n}_{k}}+
(\z^{-1}-q)P_{0}{^{ij}_{kl}})\]
belong to $\mbox{Mor}(u\cont\otimes u)$. Let denote $P_{1}:=P_{-}$,
$P_{2}:=P_{+}$ for $O_{q}(N)$ and $P_{1}:=P_{+}$, $P_{2}:=P_{-}$ for
$Sp_{q}(N)$. Then, by Lemma \ref{l-adinvsubalg},
the subspaces $\Upsilon^{0}$, $\Upsilon^{1}$ and
$\Upsilon^{2}$ of $\Gammalinv$ spanned by the sets
$\{\omega^{0}:=\frac{1}{\s}B^{\rm t}{^{m}_{k}}C^{kn}\omega_{mn}\}$,
$\{\omega^{1}_{ij}:=P_{1}(\omega_{ij})\,|\,i,j=1,\ldots,N\}$ and
$\{\omega^{2}_{ij}:=P_{2}(\omega_{ij})\,|\,i,j=1,\ldots,N\}$
respectively, are $\Delta_{R}$-invariant. Moreover, the subspaces
${\cal Y}^{0}$, ${\cal Y}^{1}$ and ${\cal Y}^{2}$ of the quantum Lie
algebra ${\cal X}$ generated by the sets $\{Y^{0}:=\frac{1}{\s}
\sum_{m}X_{mm}\}$, $\{Y^{1}_{ij}:=P_{1}{^{ij}_{kl}}X_{kl}\,|\,i,j=1,
\ldots,N\}$ and $\{Y^{2}_{ij}:=P_{2}{^{ij}_{kl}}X_{kl}\,|\,i,j=1,
\ldots,N\}$, respectively, are ad-invariant. Similarly as in case
of $SL_{q}(N)$, we obtain the following formula for the
differentiation
\[\dif a=\sum_{i,j=1}^{N}(Y^{1}_{ij}*a)\omega^{1}_{ij}+
\sum_{i,j=1}^{N}(Y^{2}_{ij}*a)\omega^{2}_{ij}+\s(Y^{0}*a)\omega^{0},
\quad a\in{\cal A}.\]

To investigate the classical limit of the structures appearing in this
article we keep the basis $\{\omega_{ij}\,|\,i,j=1,\ldots,N\}$ fixed.

Firstly let ${\cal A}=SL_{q}(N)$. We always consider the classical
limit in the sense that $z\to1$ and $q\to1$, where $z$ is the $N$-th
root of $q^{2}$. (This is not the only possibility, see
\cite{a-HSS1}.) For simplicity we shall write $\lim_{q\to1}$ for this
classical limit. Then, as shown in \cite{a-HSS1}, all functionals
$X_{ij}$, $Y^{0}$ and $Y^{1}_{ij}$ have limits when $q\to1$. It is
easily seen that the 1-forms $\omega^{0}$ and $\omega^{1}_{ij}$ and
the projections $P_{0}$ and $P_{1}$ have limits as well for $q\to1$.
We denote this limits by $\gwY^{0}$, $\gwY^{1}_{ij}$, $\gwomega^{0}$,
$\gwomega^{1}_{ij}$, $\gwP_{0}$ and $\gwP_{1}$, respectively.
It was proved in \cite{a-HSS1} that the functionals
$\gwY^{1}_{ij}$ equipped with the limit of the bracket
$[\cdot,\cdot]$ are generators of the Lie algebra $sl(N)$.

Now let ${\cal A}=O_{q}(N)$ or ${\cal A}=Sp_{q}(N)$. As proved in
\cite{a-HSS1}, for both calculi $\Gamma_{+}$ and $\Gamma_{-}$ all
functionals $X_{ij}$, $Y^{0}$, $Y^{1}_{ij}$ and
$Y^{2}_{ij}$, $i,j=1,\ldots,N$ admit limits when $q\to1$. Also, the
1-forms $\omega^{0}$, $\omega^{1}_{ij}$, $\omega^{2}_{ij}$ and the
projections $P_{k}$, $k=0,1,2$ have limits as $q\to1$. We shall use
the notations $\gwY^{0}:=\lim_{q\to1}Y^{0}$, $\gwY^{1}_{ij}:=
\lim_{q\to1}Y^{1}_{ij}$, $\gwY^{2}_{ij}:=\lim_{q\to1}Y^{2}_{ij}$,
$\gwP_{k}:=\lim_{q\to1}P_{k}$, $k=0,1,2$, $\gwomega^{0}:=
\lim_{q\to1}\omega^{0}$, $\gwomega^{1}_{ij}:=
\lim_{q\to1}\omega^{1}_{ij}$ and $\gwomega^{2}_{ij}:=
\lim_{q\to1}\omega^{2}_{ij}$ for $i,j=1,\ldots,N$. For the FODC
$\Gamma=\Gamma_{+}$ the functionals $\gwY^{1}_{ij}$, $i,j=1,\ldots,N$
equipped with the limit of the bracket $[\cdot,\cdot]$
span the Lie algebras $o(N)$ and $sp(N)$, respectively.

\section{Metrics} \label{sec-metrik}

We begin with some definitions for a general bicovariant FODC
$\tilde{\Gamma}$ over a Hopf algebra ${\cal A}$.

\begin{defin} \label{d-metrik}
A bilinear map $g:\tilde{\Gamma}\otimes_{\cal A}\tilde{\Gamma}
\rightarrow{\cal A}$ is called a {\bf metric} on $\tilde{\Gamma}$ if
$g$ is nondegenerate (i.\,e.\ $g(\xi\otimes\zeta)=0$ for all
$\zeta\in\tilde{\Gamma}$ implies $\xi=0$, $g(\xi\otimes\zeta)=0$ for
all $\xi\in\tilde{\Gamma}$ implies $\zeta=0$) and if
\begin{equation} \label{eq-alinearmetrik}
g(a\xi\otimes\zeta)=ag(\xi\otimes\zeta) \mbox{\ for any\ }a\in{\cal A},
\quad\xi\mbox{\ and\ }\zeta\in\tilde{\Gamma}.
\end{equation}
We call a metric $g$ {\bf symmetric} if $g\tilde{\sigma}=g$.
\end{defin}

For the ``ordinary'' differential calculus on Lie groups the braiding
map $\tilde{\sigma}$ is just the flip operator, so we obtain the usual
notion of a symmetric metric in this case.
{}From condition (\ref{eq-alinearmetrik}) in the preceding definition
it follows easily that a metric $g$ on $\tilde{\Gamma}$ is already
completely determined by the elements $g(\eta_{i}\otimes\eta_{j})$,
$i,j\in{\cal I}$, of ${\cal A}$.

For any metric $g$ on $\tilde{\Gamma}$ we have for all
$\xi,\zeta\in\tilde{\Gamma}$ and $a\in{\cal A}$,
\begin{equation}\label{eq-llinearmet}
\begin{gathered}
(\mbox{id}\otimes\varepsilon g)(\Delta_{L}(a\xi\otimes\zeta))=
a(\mbox{id}\otimes\varepsilon g)(\Delta_{L}(\xi\otimes\zeta)),\\
(\varepsilon g\otimes\mbox{id})(\Delta_{R}(a\xi\otimes\zeta))=
a(\varepsilon g\otimes\mbox{id})(\Delta_{R}(\xi\otimes\zeta)).
\end{gathered}
\end{equation}

\begin{defin} \label{d-invmetrik}
A metric $g$ on $\tilde{\Gamma}$ is called {\bf invariant} if for all
$\xi,\zeta\in\tilde{\Gamma}$,
\begin{equation} \label{eq-invmetrik}
(\mbox{id}\otimes\varepsilon g)(\Delta_{L}(\xi\otimes\zeta))=
g(\xi\otimes\zeta)\qquad\mbox{and}\qquad
(\varepsilon g\otimes\mbox{id})(\Delta_{R}(\xi\otimes\zeta))
=g(\xi\otimes\zeta).
\end{equation}
\end{defin}

By (\ref{eq-llinearmet}), the above definition is compatible with the
left ${\cal A}$-module structure of $\tilde{\Gamma}\otimes_{\cal A}
\tilde{\Gamma}$.

Using the representation theory of quantum groups we now show that the
invariant metrics on $\Gamma$ form a 2-parameter family for
$SL_{q}(N)$ and a 3-parameter family for $O_{q}(N)$ and $Sp_{q}(N)$.

\begin{lemma} \label{l-invmetrika}
Let ${\cal A}=SL_{q}(N)$ and let $\Gamma$ be as in Section
{\rm \ref{sec-jurcokalkule}}. A metric $g$ on $\Gamma$ is invariant if
and only if with complex parameters $\alpha$ and $\beta$
such that $\alpha\not=0,\alpha+\s\beta\not=0$
\[g(\eta_{ij}\otimes\eta_{kl})=q^{2j}\alpha\delta_{il}\delta_{jk}
+\beta\delta_{ij}\delta_{kl}.\]
\end{lemma}

\begin{bew}{}
Suppose $g$ is an invariant metric on $\Gamma$. {}From the first
equation in (\ref{eq-invmetrik})
it follows that $g(\xi\otimes\zeta)=\varepsilon g(\xi\otimes\zeta)$ for
$\xi,\zeta\in\Gammalinv$. {}From (\ref{eq-rechtskowirkung}) and
the second equation in (\ref{eq-invmetrik}) we conclude that
$g\in\mbox{Mor}(u\cont\otimes u\otimes u\cont\otimes u,1)$. Since $q$ is
not a root of unity by assumption, the multiplicities of irreducible
components in the decomposition of the tensor product representation
$u\cont\otimes u\otimes u\cont\otimes u$ are the same as in the classical
case. Therefore,
$\mbox{dim Mor}(u\cont\otimes u\otimes u\cont\otimes u,1)=2$.
Since $q^{2j}u^{j}_{n}\kappa(u^{r}_{j})=q^{2n}\delta_{nr}$ by
\cite{a-FadResTak1}, the transformations $T=(T_{ijkl})$ and
$S=(S_{ijkl})$ with
$T_{ijkl}:=q^{2j}\delta_{il}\delta_{jk}$ and
$S_{ijkl}:=\delta_{ij}\delta_{kl}$ belong to
$\mbox{Mor}(u\cont\otimes u\otimes u\cont\otimes u,1)$.
Thus we get $g(\eta_{ij}\otimes\eta_{kl})=
q^{2j}\alpha\delta_{il}\delta_{jk}+\beta\delta_{ij}\delta_{kl}$
for some complex numbers $\alpha$ and $\beta$.
This map is nondegenerate if and only if $\alpha\not=0$ and
$\alpha+\s\beta\not=0.$\\
Conversely, it is easily seen that the above formula defines an invariant
metric $g$ on $\Gamma$.
\end{bew}

\begin{lemma} \label{l-invmetrikbcd}
Let ${\cal A}=O_{q}(N)$ or ${\cal A}=Sp_{q}(N)$ and $\Gamma$ denote one
of the FODC from Section {\rm \ref{sec-jurcokalkule}}. A metric $g$ on
$\Gamma$ is invariant if and only if with complex parameters $\alpha$,
$\beta$ and $\gamma$ such that
$\alpha+\z\beta+\s\gamma\not=0$,
$\alpha-q\beta\not=0$, $\alpha+q^{-1}\beta\not=0$
\[g(\eta_{ij}\otimes\eta_{kl})=((\alpha B_{\rm 14}B_{\rm 23}+
\beta B_{\rm 12}B_{\rm 34}\Rda_{\rm 23}+\gamma B_{\rm 12}B_{\rm 34})
C{^{\rm t}_{\rm 1}}C{^{\rm t}_{\rm 3}})_{ijkl}.\]
\end{lemma}

In Lemma \ref{l-invmetrikbcd} we used the following notation. Let
$C=(C^{i}_{j})$ be the matrix of the metric which occurs in the
definition of $O_{q}(N)$ resp.\ $Sp_{q}(N)$ (see \cite{a-FadResTak1})
and let $B=(B^{i}_{j})$ its inverse matrix. Set $C^{\rm t}:=
(C^{\rm t}{^{i}_{j}})=(C^{ji})$\,$(=(C^{j}_{i}))$ and
$B^{\rm t}:=(B^{\rm t}{^{i}_{j}})=(B_{ji})$\,$(=(B^{j}_{i}))$.
Then the notation in Lemma \ref{l-invmetrikbcd} is the
usual leg numbering notation, i.\,e.\ the equation therein reads as
\[g(\eta_{ij}\otimes\eta_{kl})=(\alpha B_{ml}B_{jn}+
\beta B_{mr}B_{sl}\Rda{^{rs}_{jn}}+\gamma B_{mj}B_{nl})
C^{\rm t}{^{m}_{i}}C^{\rm t}{^{n}_{k}}.\]

\begin{bew}{of Lemma \ref{l-invmetrikbcd}}
The proof is similar to the proof of Lemma \ref{l-invmetrika}.
The decomposition of the tensor product gives now
$\mbox{dim Mor}(u\cont\otimes u\otimes u\cont\otimes u,1)=3$.
The conditions on the coefficients ensure the nondegeneracy of
the metric $g$.
\end{bew}

Some straightforward computations based on Lemma \ref{l-sigmamatrix} and
the particular form of the invariant metrics in Lemma \ref{l-invmetrika}
and \ref{l-invmetrikbcd} prove the following

\begin{folg} \label{f-symmmetrik}
Any invariant metric on $\Gamma$ is symmetric.
\end{folg}

Let us consider again a general bicovariant FODC $\tilde{\Gamma}$ over
an arbitrary Hopf algebra ${\cal A}$.

Let $g$ be a metric on $\tilde{\Gamma}\otimes_{\cal A}\tilde{\Gamma}$
and let $g_{ij}$ be the matrix coefficients of $g$, i.\,e.\
$g(\eta_{i}\otimes\eta_{j})=g_{ij}$ with respect to a fixed basis
$\{\eta_{i}\}$ of $\linv{\tilde{\Gamma}}$. Recall that $\{\chi_{i}\}$
is the dual basis of $\{\eta_{i}\}$. We suppose that $\{a_{j}\,|\,j
\in{\cal I}\}$ is a finite subset of ${\cal A}$ such that $\chi_{i}(a_{j})=
\delta_{ij}$.

\begin{defin}
A map $g^{*}:\tilde{\cal X}\otimes\tilde{\cal X}\rightarrow{\cal A}$
such that $g^{*}(\chi_{i}\otimes\chi_{j})=g^{*}{^{ij}}$ and
$\sum_{k}g_{ik}g^{*}{^{kj}}=\delta{^{j}_{i}}$
for all $i,j\in{\cal I}$ is called the {\bf dual metric of} $g$.
The metric $g^{*}$ is called {\bf ad-invariant} if
\begin{equation} \label{eq-adinvmetrik}
g^{*}([\chi',f_{(1)}]\otimes[\chi'',f_{(2)}])=
f(1)g^{*}(\chi'\otimes\chi'')\quad
\text{for all $\chi',\chi''\in\tilde{\cal X}$, $f\in{\cal A}^{0}$.}
\end{equation}
\end{defin}

Note that the nondegeneracy of $g$ corresponds to the nondegeneracy
of $g^{*}$.

There is an interesting link between invariant and ad-invariant
metrics given by

\begin{satz} \label{s-dualmetrikinvarianz}
If $g$ is an invariant metric on $\tilde{\Gamma}$ then the dual
metric $g^{*}$ of $g$ is ad-invariant.\\
Conversely, let $g^{*}:\tilde{\cal X}\otimes\tilde{\cal X}
\to\mathbb{C}$ be an ad-invariant metric and suppose
that ${\cal A}^{0}$ separates the points of ${\cal A}$.
Then $g$ is an invariant metric.
\end{satz}

\begin{bew}{}
In the proof we take the bases of $\linv{\tilde{\Gamma}}$ and
$\tilde{\cal X}$ described above.

Let $g$ be invariant. {}From the first equation in
(\ref{eq-invmetrik}) it follows that
$g(\eta_{i}\otimes\eta_{j})=g_{ij}\in\mathbb{C}$ for all $i,j\in{\cal I}$.
Since $g$ is nondegenerate, the dual metric $g^{*}$ is defined and is
a map to $\mathbb{C}\cdot 1$. The second equation in
(\ref{eq-invmetrik}) means that
$g_{ij}\tilde{v}^{i}_{m}\tilde{v}{^{j}_{n}}=g_{mn}$ which yields
$g^{*}{^{ij}}\tilde{v}{^{m}_{i}}\tilde{v}{^{n}_{j}}=g^{*}{^{mn}}$ for
the dual metric $g^{*}$. Applying functionals $f\in{\cal A}^{0}$ to the
last equation we get $g^{*}(f_{(1)}(\tilde{v}^{m}_{i})\chi_{i}\otimes
f_{(2)}(\tilde{v}^{n}_{j})\chi_{j})=f(1)g^{*}(\chi_{m}\otimes
\chi_{n})$. Using
$[\chi_{i},f]=f(\tilde{v}^{i}_{k})\chi_{k}$ for all $f\in{\cal A}^0$
the latter is equivalent to
$g^{*}([\chi_{m},f_{(1)}]\otimes[\chi_{n},f_{(2)}])=f(1)g^{*}(\chi_{m}
\otimes\chi_{n})$ for all $m,n\in{\cal I}$ and all $f\in{\cal A}^{0}$,
i.\,e.\ the dual metric is ad-invariant.

Suppose that $g^{*}:\tilde{\cal X}\otimes\tilde{\cal X}\to\mathbb{C}$
is ad-invariant. Then the matrix elements of $g$ are also
complex numbers. This implies the first equation in
(\ref{eq-invmetrik}). Reversing the
reasoning from the preceding paragraph, it follows from the
ad-invariance of $g^{*}$ that $f(g^{*}{^{ij}}\tilde{v}{^{m}_{i}}
\tilde{v}{^{n}_{j}}-g^{*}{^{mn}})=0$ for all $m,n\in{\cal I}$ and $f\in
{\cal A}^{0}$. Since ${\cal A}^{0}$ separates the points of
${\cal A}$, we obtain $g^{*}{^{ij}}\tilde{v}{^{m}_{i}}
\tilde{v}{^{n}_{j}}=g^{*}{^{mn}}$ for all $m,n\in{\cal I}$ from which the
second equation in (\ref{eq-invmetrik}) follows.
This proves the invariance of $g$.
\end{bew}

Next we specialize again to the bicovariant FODC $\Gamma=
\Gamma_{\pm,z}$ over ${\cal A}=SL_{q}(N)$, $O_{q}(N)$, $Sp_{q}(N)$.
We compute the dual metrics $g^{*}$ of the invariant metrics $g$
from Lemmas \ref{l-invmetrika} and \ref{l-invmetrikbcd} and their
restrictions to the ad-invariant
subspaces of the quantum Lie algebra ${\cal X}$. Let us say that two
subspaces $\Upsilon$ and $\Upsilon'$ of $\Gammalinv$ are {\bf mutually
orthogonal with respect to a metric} $g$ on $\Gammalinv$ if $g(x
\otimes x')=g(x'\otimes x)=0$ for all $x\in\Upsilon$,
$x'\in\Upsilon'$. A similar notion is used for the dual metric $g^{*}$
on ${\cal X}$.

For $SL_{q}(N)$ we have (see Lemma \ref{l-invmetrika})
$g(\eta_{ij}\otimes\eta_{kl})=q^{2j}\alpha\delta_{il}\delta_{jk}+
\beta\delta_{ij}\delta_{kl}$ and we get
\begin{equation} \label{eq-dualmeta}
g^{*}(\chi_{ij}\otimes\chi_{kl})=\left(q^{-2i}\alpha^{-1}
\delta_{il}\delta_{jk}-\frac{\beta}{\alpha(\alpha+\s\beta)}q^{-2i-2k}
\delta_{ij}\delta_{kl}\right)
\end{equation}
for all $i,j,k,l=1,\ldots,N$.
For $O_{q}(N)$ and $Sp_{q}(N)$ the dual metric of the metric in Lemma
\ref{l-invmetrikbcd} is
\[g^{*}(\chi_{ij}\otimes\chi_{kl})=\frac{B^{\rm t}{^{i}_{m}}
B^{\rm t}{^{k}_{n}}}{(\alpha\!-\!q\beta)(\alpha\!+\!q^{-1}\beta)}
\left(\alpha C_{ml}C_{jn}\!-\!\beta\Rdam{^{jn}_{rs}}C_{mr}C_{sl}+
\frac{-Q\alpha\beta\!-\!\alpha\gamma\!+\!\z^{-1}\beta\gamma}
{\alpha\!+\!\z\beta\!+\!\s\gamma}C_{mj}C_{nl}\right)\]
for $i,j,k,l=1,\ldots,N$.

{}From equations (\ref{eq-chifunkt}) and (\ref{eq-lfunktionale}) we
get the transformation formula between the bases $\{X_{ij}\}$ and
$\{\chi_{ij}\}$. Then we can express the generators of the
ad-invariant subspaces in terms of the basis $\{\chi_{ij}\}$.
For the quantum group $SL_{q}(N)$ and the FODC
$\Gamma=\Gamma_{\pm,z}$ we obtain the formulas
$Y^{0}=\s^{-1}\mu_{\pm,z}^{-1}\sum_{k}\chi_{kk}$ and
$Y^{1}_{ij}=\nu_{\pm,z}^{-1}P_{1}{^{ij}_{kl}}\chi_{kl}$. For the
quantum groups $O_{q}(N)$ and $Sp_{q}(N)$ and the calculus
$\Gamma=\Gamma_{\pm1}$ we get
$Y^{0}=\s^{-1}\mu_{0}{_{\pm}}^{-1}\sum_{m}\chi_{mm}$,
$Y^{1}_{ij}=\mu_{1}{_{\pm}}^{-1}P_{1}{^{ij}_{kl}}\chi_{kl}$ and
$Y^{2}_{ij}=\mu_{2}{_{\pm}}^{-1}P_{2}{^{ij}_{kl}}\chi_{kl}$. Here
the constants are defined by
$\mu_{+,z}=\s(z^{-1}-1)+z^{-1}q^{-1}Q$,
$\mu_{-,z}=\s(z-1)-zq^{-2N-1}Q$,
$\nu_{+,z}=z^{-1}q^{-1}Q$,
$\nu_{-,z}=-zq^{-2N-1}Q$,
$\mu_{0}{_{+}}=(\z-\z^{-1})Q$,
$\mu_{0}{_{-}}=(\z-\z^{-1})Q+2\s$,
$\mu_{1}{_{\pm}}=\pm\z(1+q^{-N})Q$ and
$\mu_{2}{_{\pm}}=\pm(\z-\z^{-1}q^{N})Q$.

One easily verifies that ${\cal Y}^{0}$ and ${\cal Y}^{1}$ are
orthogonal subspaces of ${\cal X}$ with respect to the dual metric
$g^{*}$. The restrictions of $g^{*}$ to ${\cal Y}^{0}\otimes
{\cal Y}^{0}$ and ${\cal Y}^{1}\otimes{\cal Y}^{1}$ are given by
\[g^{*}(Y^{0}\otimes Y^{0})=\frac{1}{\s(\alpha+\s\beta)
\mu_{\pm,z}^{2}},\qquad g^{*}(Y^{1}_{ij}\otimes Y^{1}_{kl})=
\frac{1}{\alpha\nu_{\pm,z}^{2}}\left(q^{-2i}\delta_{il}\delta_{jk}-
\frac{q^{-2i-2k}}{\s}\delta_{ij}\delta_{kl}\right)\]
for all $i,j,k,l=1,\ldots,N$.

Then $\Upsilon^{0}$ and $\Upsilon^{1}$ are orthogonal subspaces of
$\Gammalinv$ with respect to all invariant metrics $g$. If
$g$ has the form as in Lemma \ref{l-invmetrika}, then we have
\[g(\omega^{0}\otimes\omega^{0})=\frac{(\alpha+\s\beta)
\mu_{\pm,z}^{2}}{\s},\qquad g(\omega^{1}_{ij}\otimes\omega^{1}_{kl})=
\alpha\nu_{\pm,z}^{2}\left(q^{2j}\delta_{il}\delta_{jk}-\frac{1}{\s}
\delta_{ij}\delta_{kl}\right)\]
for all $i,j,k,l=1,\ldots,N$.

The ad-invariant subspaces ${\cal Y}^{0}$, ${\cal Y}^{1}$ and
${\cal Y}^{2}$ of the quantum Lie algebras ${\cal X}$ of $O_{q}(N)$
and $Sp_{q}(N)$ are mutually orthogonal with respect to $g^{*}$. The
restrictions of $g^{*}$ to ${\cal Y}^{i}\otimes{\cal Y}^{i}$,
$i=0,1,2$ are described by the formulas
\[g^{*}(Y^{0}\otimes Y^{0})=
\frac{1}{\mu_{0}{_{\pm}}^{2}\s\alpha_{0}},\]
\[g^{*}(Y^{1}_{ij}\otimes Y^{1}_{kl})=\frac{B^{\rm t}{^{i}_{m}}
B^{\rm t}{^{k}_{n}}}{\mu_{1}{_{\pm}}^{2}(\z^{-1}q^{N}\!+\!\z q^{-N})
\alpha_{1}}\left(\z^{-1}q^{N}C^{ml}C^{jn}\!-\!
\Rdam{^{jn}_{rs}}C^{mr}C^{sl}\!+\!\frac{\z^{-1}(1\!-\!q^{N})}{\s}
C^{mj}C^{nl}\right),\]
\[g^{*}(Y^{2}_{ij}\otimes Y^{2}_{kl})=\frac{B^{\rm t}{^{i}_{m}}
B^{\rm t}{^{k}_{n}}}{\mu_{2}{_{\pm}}^{2}(\z^{-1}q^{N}\!+\!\z q^{-N})
\alpha_{2}}\left(\z q^{-N}C^{ml}C^{jn}\!-\!\Rdam{^{jn}_{rs}}C^{mr}C^{sl}
\!-\!\frac{\z^{-1}\!+\!\z q^{-N}}{\s}C^{mj}C^{nl}\right),\]
where we use the abbreviations
$\alpha_{0}:=\alpha+\z\beta+\s\gamma$,
$\alpha_{1}:=\alpha-\z^{-1}q^{N}\beta$,
$\alpha_{2}:=\alpha+\z q^{-N}\beta$.
Note that all denominators in the above formulas are
non-zero, since the metric is nondegenerate and $\Gamma$ is
a FODC as in \cite{a-SchSch1}.

The corresponding subspaces $\Upsilon^{0}$, $\Upsilon^{1}$ and
$\Upsilon^{2}$ of $\Gammalinv$ are mutually orthogonal with respect to
all invariant metrics $g$. If $g$ is as in Lemma \ref{l-invmetrikbcd},
then we have
\[g(\omega^{0}\otimes\omega^{0})=\frac{\mu_{0}{_{\pm}}^{2}\alpha_{0}}
{\s},\]
\[g(\omega^{1}_{ij}\otimes\omega^{1}_{kl})=
\frac{\mu_{1}{_{\pm}}^{2}\alpha_{1}}{\z^{-1}q^{N}+\z q^{-N}}\left(
\z q^{-N}B_{ml}B_{jn}-B_{mr}B_{sl}\Rda{^{rs}_{jn}}+\frac{\z(1-q^{-N})}
{\s}B_{mj}B_{nl}\right)C^{\rm t}{^{m}_{i}}C^{\rm t}{^{n}_{k}},\]
\[g(\omega^{2}_{ij}\otimes\omega^{2}_{kl})=
\frac{\mu_{2}{_{\pm}}^{2}\alpha_{2}}{\z^{-1}q^{N}+\z q^{-N}}\left(
\z^{-1}q^{N}B_{ml}B_{jn}+B_{mr}B_{sl}\Rda{^{rs}_{jn}}-
\frac{\z+\z^{-1}q^{N}}{\s}B_{mj}B_{nl}\right)C^{\rm t}{^{m}_{i}}
C^{\rm t}{^{n}_{k}}\]
for all $i,j,k,l=1,\ldots,N$.

Now we want to examine the classical limits.

Let ${\cal A}=SL_{q}(N)$ and let $g$ be an invariant metric as
described in Lemma \ref{l-invmetrika}. The complex numbers
$\alpha$ and $\beta$ may of course depend
on the parameters $q$ and $z$.
Let us assume that the functions $\alpha=\alpha(q,z)$ and
$\beta=\beta(q,z)$ are choosen such that the limits
$c_{1}:=\lim_{q\to1}Q^{2}\alpha(q,z)$ and
$c_{0}:=\lim_{q\to1}Q^{4}(\alpha(q,z)+\s\beta(q,z))$ exist and are
non-zero. Then it follows immediately from the existence of the
classical limits as discussed at the end of Section
\ref{sec-jurcokalkule} that the invariant metric $g$ and its dual
metric admit limits when $q\to1$ and $z\to1$. The restriction
of the limit of $g^{*}$ to the linear functionals $\gwY^{1}_{ij}$,
$i,j=1,\ldots,N$, is just a complex multiple of the Killing form for
$sl(N)$. Similar results are valid for the quantum groups $O_{q}(N)$
and $Sp_{q}(N)$ and for the FODC $\Gamma_{+}$ if we suppose that
for the functionals $\alpha_{k}=\alpha_{k}(q)$, $k=0,1,2$ the
limits $c_{0}:=\lim_{q\to1}\mu_{0}{_{+}}^{2}\alpha_{0}(q)$,
$c_{1}:=\lim_{q\to1}\mu_{1}{_{+}}^{2}\alpha_{1}(q)$
and $c_{2}:=\lim_{q\to1}\mu_{2}{_{+}}^{2}\alpha_{2}(q)$
exist and are non-vanishing.

\section{Connections} \label{sec-zusammenhang}

We begin with some general definitions (cf.\ \cite{b-Connes1})
which, of course, apply to any differential calculus
$\tilde{\Gamma}^{\wedge}=\bigoplus_{n=0}^{\infty}\tilde{\Gamma}^{n}$
over an arbitrary algebra ${\cal A}$.

Let ${\cal E}$ be a left ${\cal A}$-module. A {\bf left connection} on
${\cal E}$ is a linear map
$\nabla:{\cal E}\rightarrow\tilde{\Gamma}\otimes_{\cal A}{\cal E}$
such that
\begin{equation} \label{eq-deflzush}
\nabla(a\zeta)=\dif aPotimes\zeta+a\nabla(\zeta)\qquad
\mbox{for all }a\in{\cal A}\mbox{ and }\zeta\in{\cal E}.
\end{equation}
Let $\tilde{\Gamma}{\cal E}:=\tilde{\Gamma}^{\wedge}\otimes_{\cal A}
{\cal E}$ be the ``${\cal E}$-valued differential forms''. A
connection $\nabla$ on ${\cal E}$ admits a unique extension to a
linear map $\nabla:\tilde{\Gamma}{\cal E}\rightarrow\tilde{\Gamma}
{\cal E}$ of degree one such that
\[\nabla(\alpha\zeta)=(\dif\alpha)\zeta+(-1)^{n}\alpha\nabla(\zeta)
\qquad\mbox{for }\alpha\in\tilde{\Gamma}^{n}\mbox{\ and\ }\zeta\in
\tilde{\Gamma}{\cal E}.\]
The mapping $R(\nabla)=\nabla^{2}:{\cal E}\rightarrow
\tilde{\Gamma}^{2}\otimes_{\cal A}{\cal E}$ is called the
{\bf curvature} of the connection $\nabla$.
Clearly, $R(\nabla)$ is ${\cal A}$-linear, i.\,e.\ $R(\nabla)(a\zeta)=
aR(\nabla)(\zeta)$ for $a\in{\cal A}$ and $\zeta\in{\cal E}$.

Similar concepts can be defined for a right ${\cal A}$-module ${\cal E}$.
A right connection on ${\cal E}$ is then a linear map $\nabla:{\cal E}
\rightarrow{\cal E}\otimes_{\cal A}\tilde{\Gamma}$ satisfying
$\nabla(\zeta a)=\zeta\otimes\dif a+\nabla(\zeta)a$ for $a\in{\cal A}$
and $\zeta\in{\cal E}$. It extends uniquely to a linear map
$\nabla:{\cal E}\otimes_{\cal A}\tilde{\Gamma}^{\wedge}\rightarrow
{\cal E}\otimes_{\cal A}\tilde{\Gamma}^{\wedge}$ such that
$\nabla(\zeta\alpha)=(-1)^{n}\zeta\dif\alpha+\nabla(\zeta)\alpha$,
$\alpha\in\tilde{\Gamma}^{\wedge}$ and $\zeta\in{\cal E}
\otimes_{\cal A}\tilde{\Gamma}^{n}.$ The curvature of $\nabla$ is
$R(\nabla):=\nabla^{2}:{\cal E}\rightarrow{\cal E}\otimes_{\cal A}
\tilde{\Gamma}^{2}$.

If $\nabla$ is a connection on a left ${\cal A}$-module ${\cal E}$, then
there is a connection $\nabla^{*}$ on the right ${\cal A}$-module
${\cal E}^{*}$, called the {\bf dual connection} of $\nabla$, defined by
\[<\xi,\nabla^{*}(\zeta)>=\dif<\xi,\zeta>-<\nabla(\xi),\zeta>\qquad
\mbox{for }\xi\in{\cal E},\quad\zeta\in{\cal E}^{*}.\]
By a {\bf right} (resp.\ {\bf left}) {\bf connection} on an
${\cal A}$-bimodule ${\cal E}$ we mean a connection on ${\cal E}$ when
${\cal E}$ is considered as a right (resp.\ left) ${\cal A}$-module.

We now specialize to the case of our main interest where
${\cal E}=\tilde{\Gamma}$ is considered as a left ${\cal A}$-module.
Suppose that $\nabla:\tilde{\Gamma}\rightarrow\tilde{\Gamma}
\otimes_{\cal A}\tilde{\Gamma}$ is a (left) connection on
$\tilde{\Gamma}$.
Then the {\bf torsion} $T(\nabla)$ is defined by $T(\nabla)=\dif-
\mbox{m}\nabla$, where $\mbox{m}:\tilde{\Gamma}\otimes_{\cal A}
\tilde{\Gamma}\rightarrow\tilde{\Gamma}^{2}$ denotes the
multiplication map, i.\,e.\ $\mbox{m}(\xi\otimes\zeta)=
\xi\wedge\zeta$ for $\xi,\zeta\in\tilde{\Gamma}$. The torsion
$T(\nabla)$ is ${\cal A}$-linear, since
$T(\nabla)(a\xi)=\dif(a\xi)-\mbox{m}\nabla(a\xi)=
\dif a\wedge\xi+a\dif\xi-da\wedge\xi-a\mbox{m}\nabla(\xi)=
aT(\nabla)(\xi)$ for $a\in{\cal A},\,\xi\in\tilde{\Gamma}$ by the
Leibniz rule and (\ref{eq-deflzush}).

\begin{defin} \label{d-bikovzush}
The connection $\nabla$ is called {\bf bicovariant} if
\begin{equation} \label{eq-invzush}
\Delta_{L}\nabla=(\mbox{id}\otimes\nabla)\Delta_{L}\qquad\mbox{and}
\qquad\Delta_{R}\nabla=(\nabla\otimes\mbox{id})\Delta_{R}.
\end{equation}
\end{defin}

A simple computation shows that the set of all bicovariant connections
on $\tilde{\Gamma}$ forms a complex affine space $\mbox{BC}(
\tilde{\Gamma})$.

Let $\tilde{\Gamma}^{\wedge}=\bigoplus_{n=0}^{\infty}
\tilde{\Gamma}^{n}$be a bicovariant differential calculus over a Hopf
algebra ${\cal A}$. Let $\{\eta_{i}\,|\,i\in{\cal I}\}$ be a basis in
$\linv{\tilde{\Gamma}}$, $\Delta_{R}(\eta_{i})=\eta_{j}\otimes
\tilde{v}^{j}_{i}$, and let $D(\nabla){_{i}^{jk}}$ denote arbitrary
elements of ${\cal A}$. Then the map $\nabla:\linv{\tilde{\Gamma}}
\rightarrow\tilde{\Gamma}\otimes_{\cal A}\tilde{\Gamma}$, defined by
$\nabla(\eta_{i})=D(\nabla){_{i}^{jk}}\eta_{j}\otimes\eta_{k}$ extends
uniquely to a connection on $\tilde{\Gamma}$ and any connection on
$\tilde{\Gamma}$ is of this form.

\begin{lemma} \label{l-Dmorph}
Let $\tilde{v}$ denote the corepresentation of ${\cal A}$ defined by
$\Delta_{R}$ on $\linv{\tilde{\Gamma}}$.
The connection $\nabla$ on $\tilde{\Gamma}$ is bicovariant if and only if
$D(\nabla){_{i}^{jk}}\in\mathbb{C}$ for all $i,j,k$ and
$D(\nabla)=(D(\nabla)_{i}^{jk})\in\mbox{Mor}(\tilde{v},\tilde{v}\otimes
\tilde{v})$.
\end{lemma}

\begin{bew}{}
Let $\nabla$ be a bicovariant connection on $\tilde{\Gamma}$ and
$\nabla(\eta_{i})=D(\nabla){_{i}^{jk}}\eta_{j}\otimes\eta_{k}$,
$D(\nabla){_{i}^{jk}}\in{\cal A}$. Because of the first formula of
(\ref{eq-invzush}) we
have $\Delta(D(\nabla){_{i}^{jk}})=1\otimes D(\nabla){_{i}^{jk}}$ and
so $D(\nabla){_{i}^{jk}}\in\mathbb{C}\cdot 1$.
The second equation of (\ref{eq-invzush}) tells us that
$(\nabla\otimes\mbox{id})\Delta_{R}(\eta_{i})=
\eta_{j}\otimes\eta_{k}\otimes D(\nabla){_{x}^{jk}}\tilde{v}{^{x}_{i}}=
\Delta_{R}\nabla(\eta_{i})=\eta_{j}\otimes\eta_{k}\otimes
(\tilde{v}\otimes\tilde{v}){^{jk}_{xy}}D(\nabla){_{i}^{xy}}$, i.\,e.\
$D(\nabla)\tilde{v}=(\tilde{v}\otimes\tilde{v})D(\nabla)$.

Let now $D(\nabla)_{i}^{jk}\in\mathbb{C}$ and $D(\nabla)\in\mbox{Mor}
(\tilde{v},\tilde{v}\otimes\tilde{v})$. Then the above equations
written in reversed order show that (\ref{eq-invzush})
is true for the connection $\nabla$ defined by
$\nabla(\eta_{i})=D(\nabla){_{i}^{jk}}\eta_{j}\otimes\eta_{k}$.
\end{bew}

Now let $g$ be an invariant metric on $\tilde{\Gamma}$, cf.\
Definition \ref{d-invmetrik}. Let us introduce some new concepts.

\begin{defin} \label{d-vertrmet}
We say that the connection $\nabla$ on $\tilde{\Gamma}$ is
{\bf compatible with the metric} $g$ if
\[(\mbox{id}\otimes g)(\nabla(\xi)\otimes\zeta)+
(g\otimes\mbox{id})(\mbox{id}\otimes\tilde{\sigma})(\xi\otimes\nabla
(\zeta))-\dif g(\xi\otimes\zeta)=0\quad\mbox{for any }\xi\in
\tilde{\Gamma},\zeta\in\linv{\tilde{\Gamma}}.\]
\end{defin}

\begin{lemma} \label{l-vertrneu}
Let $g$ denote an invariant metric. Then a connection $\nabla$ on
$\tilde{\Gamma}$ is compatible with $g$ if and only if
\begin{equation} \label{eq-vertrneu}
(\mbox{id}\otimes g)(\nabla(\eta_{i})\otimes\eta_{j})+
(g\otimes\mbox{id})(\eta_{i}\otimes\tilde{\sigma}\nabla(\eta_{j}))=0
\mbox{\ for all\ }i,j\in{\cal I}.
\end{equation}
\end{lemma}

\begin{bew}{}
In both directions we use only that the metric $g$ is left-invariant,
i.\,e.\ the first equation in (\ref{eq-invmetrik}) is fulfilled.
This is equivalent to
$g_{ij}=g(\eta_{i}\otimes\eta_{j})\in\mathbb{C}$.

Suppose $\nabla$ is compatible with $g$. Because of the first
formula in (\ref{eq-invmetrik})
and Definition \ref{d-vertrmet} equation (\ref{eq-vertrneu}) is valid.
Let us now assume equation (\ref{eq-vertrneu}) is fulfilled and introduce
arbitrary elements $\xi\otimes\zeta=a_{ij}\eta_{i}\otimes\eta_{j}
\in\tilde{\Gamma}\otimes_{\cal A}\tilde{\Gamma}$, $a_{ij}\in{\cal A}$.
Then the assertion follows from the computation
\begin{eqnarray*}
\lefteqn{(\mbox{id}\otimes g)(\nabla(a_{ij}\eta_{i})\otimes\eta_{j})+
(g\otimes\mbox{id})(a_{ij}\eta_{i}\otimes\tilde{\sigma}\nabla
(\eta_{j}))-\dif g(a_{ij}\eta_{i}\otimes\eta_{j})=}\\
&& =(\mbox{id}\otimes g)((\dif a_{ij}\otimes\eta_{i}+
a_{ij}\nabla\eta_{i})\otimes\eta_{j})+
a_{ij}(g\otimes\mbox{id})(\eta_{i}\otimes\tilde{\sigma}\nabla\eta_{j})
-\\
&& -(\dif (a_{ij})g(\eta_{i}\otimes\eta_{j})+
a_{ij}\dif g(\eta_{i}\otimes\eta_{j}))=\\
&& =a_{ij}((\mbox{id}\otimes g)(\nabla\eta_{i}\otimes\eta_{j})+
(g\otimes\mbox{id})(\eta_{i}\otimes\tilde{\sigma}\nabla\eta_{j}))-
a_{ij}\dif g(\eta_{i}\otimes\eta_{j})
\end{eqnarray*}
and the assumption $g_{ij}\in\mathbb{C}$.
\end{bew}

\begin{defin} \label{d-LCZ}
Let $g$ be an invariant metric. A bicovariant connection $\nabla$ on
$\tilde{\Gamma}$ is called a {\bf Levi-Civita connection} (with
respect to the metric $g$) if $\nabla$ is compatible with $g$ and has
vanishing torsion $T(\nabla)$.
\end{defin}

Recall that any basis of the vector space $\linv{\tilde{\Gamma}}$ is a
free left ${\cal A}$-module basis of $\tilde{\Gamma}$. Therefore, any
connection on $\tilde{\Gamma}$, its curvature and its torsion are
uniquely determined by its values on such a basis.

Now we consider the pairing between $\tilde{\Gamma}$ and
$\tilde{\cal X}\otimes{\cal A}$ as a right ${\cal A}$-module
(see formula (\ref{eq-dualitaet})).
Since for a connection $\nabla$ on $\tilde{\Gamma}$ there is a dual
connection $\nabla^{*}$ on $\tilde{\cal X}\otimes{\cal A}$, it is
natural to ask what the conditions in Definition \ref{d-LCZ} mean in
terms of $\nabla^{*}$. For this we assume that $\tilde{\Gamma}^{2}=
(\tilde{\Gamma}\otimes_{\cal A}\tilde{\Gamma})/
\ker(\tilde{\sigma}-\mbox{id})$.

Suppose that $\nabla$ is a connection on $\tilde{\Gamma}$ and let
$\nabla^{*}$ be its dual connection.
For arbitrary $\xi,\zeta\in\tilde{\cal X}$ we define
$\nabla^{*}_{\xi}(\zeta)=<\nabla^{*}(\zeta),\xi>$ and mean the pairing
between $\tilde{\cal X}\otimes\tilde{\Gamma}$ and $\tilde{\cal X}$
(cf.\ (\ref{eq-dualitaet})).
For the basis elements $\eta_{i}$ of $\linv{\tilde{\Gamma}}$ we write
$\nabla(\eta_{i})=\Gamma_{i}^{jk}\eta_{j}\otimes\eta_{k}$ with
$\Gamma_{i}^{jk}\in{\cal A}$.
Recall that by Lemma \ref{l-Dmorph}, we have $\Gamma_{i}^{jk}\in
\mathbb{C}$ for all $i,j,k\in{\cal I}$.
Then the pairing gives $\nabla^{*}_{\chi_{i}}(\chi_{j})=
-\Gamma_{k}^{ij}\chi_{k}$ for the elements of the dual basis.
Let $\eta$ denote the left- and right-invariant
element of $\tilde{\Gamma}$ used in the definition of the
differentiation, i.\,e.\ $\dif a=\eta a-a\eta$ for $a\in{\cal A}$.
Recall that the torsion is ${\cal A}$-linear. Therefore, because of
$\dif\eta_{i}=\eta\wedge\eta_{i}+\eta_{i}\wedge\eta$, the torsion of
$\nabla$ is vanishing if and only if $\dif\eta_{i}=
\mbox{m}\nabla(\eta_{i})$ or equivalently
\[(\mbox{id}-\tilde{\sigma})(\eta\otimes\eta_{i}+\eta_{i}\otimes\eta-
\nabla(\eta_{i}))=0\qquad\mbox{for all\ }i\in{\cal I}.\]
Using the notation above the latter is equivalent to the equation
\[\nabla^{*}_{\chi_{i}}(\chi_{j})-
\tilde{\sigma}^{ij}_{kl}\nabla^{*}_{\chi_{k}}(\chi_{l})
=[\chi_{i},\chi_{j}]\mbox{\ for all\ }i,j\in{\cal I}.\]
Suppose that $g$ is an invariant metric. By dualizing the condition in
Lemma \ref{l-vertrneu}, it follows that $\nabla$ is compatible with
the metric $g$ if and only if
\[g^{*}(\chi_{i}\otimes\nabla^{*}_{\chi_{j}}(\chi_{k}))+g^{*}(
\tilde{\sigma}{^{ij}_{mn}}\nabla^{*}_{\chi_{m}}(\chi_{n})\otimes
\chi_{k})=0\mbox{\ for all\ }i,j,k\in{\cal I}.\]
The above equations show that our concepts are analogous to the
corresponding notions in classical differential geometry.

\section{Levi-Civita connections on $SL_{q}(N)$}
\label{sec-lcza}

In this section we examine the differential calculi $\Gamma_{\pm,z}$ for
${\cal A}=SL_{q}(N)$. After a short lemma we will prove our main results.

\begin{lemma} \label{l-zusa}
We have $\rm{dim\ BC}(\Gamma)=5$ for $N=2$ and $\rm{dim\ BC}(\Gamma)=6$
for $N\geq3$.
\end{lemma}

\begin{bew}{}
Since $q$ is not a root of unity, decompositions of tensor product
representations of ${\cal A}=SL_{q}(N)$ can be labelled by Young
tableaus as in the classical case. Therefore, we obtain
$u\cont\otimes u=[0]\oplus[2,1^{N-2}]$ and
$u\cont\otimes u\otimes u\cont\otimes u=2[0]\oplus k[2,1^{N-2}]\oplus$
other terms with $k=3$ for $N=2$ and $k=4$ for $N\geq3$. Then by the
general representation theory we conclude
that $\mbox{dim Mor}(v,v\otimes v)=5$ for $N=2$ and
$\mbox{dim Mor}(v,v\otimes v)=6$ for $N\geq3$. By Lemma \ref{l-Dmorph},
a connection $\nabla$ on $\Gamma_{\pm,z}$ is bicovariant if and only
if $D(\nabla)\in\mbox{Mor}(v,v\otimes v)$, where $v=u\cont\otimes u$.
Thus the assertion of the lemma follows.
\end{bew}

Let $\Gamma^{\wedge}=\bigoplus_{n=0}^{\infty}\Gamma^{n}$ be a differential
calculus over ${\cal A}$ which contains $\Gamma_{\pm,z}$ as its first
order differential calculus $\Gamma$. As in \cite{a-Woro2} we suppose that
$\Gamma^{2}=(\Gamma\otimes_{\cal A}\Gamma)/\ker(\sigma-\mbox{id})$. Then
we have the following

\begin{thm} \label{t-lcza}
For any invariant metric $g$ on $\Gamma$ there exists precisely one
Levi-Civita connection $\nabla$. If
$g(\eta_{ij}\otimes\eta_{kl})=q^{2j}\alpha\delta_{il}\delta_{jk}
+\beta\delta_{ij}\delta_{kl}$ with $\alpha,\beta\in\mathbb{C},
\alpha\not=0,\,\alpha+\s\beta\not=0$ then for
$\Gamma_{+,z}$ this Levi-Civita connection is given by
\begin{eqnarray}
\lefteqn{ \nabla(\eta_{ij})=\frac{Q}{2}(q^{2p-2a-2c}\Rda{^{bn}_{am}}
\Rda{^{dp}_{cn}}\Rda{^{im}_{jp}}\eta_{ab}\otimes\eta_{cd}
-q^{-1}\eta_{ia}\otimes\eta_{aj} -} \nonumber\\
& &-Q\eta_{ij}\otimes\eta-Q\eta\otimes\eta_{ij}
+\delta_{ij}Q(1+\s\alpha^{-1}\beta)q^{-2a}\eta_{ab}\otimes
\eta_{ba}-\nonumber\\
& &-\delta_{ij}(q^{2N+1}Q^{2}+Q\alpha^{-1}\beta)\eta\otimes
\eta)\nonumber
\end{eqnarray}
and for $\Gamma_{-,z}$ by
\begin{eqnarray}
\lefteqn{ \nabla(\eta_{ij})=\frac{Q}{2}(-q^{2p-2a-2c}\Rda{^{bn}_{am}}
\Rda{^{dp}_{cn}}\Rda{^{im}_{jp}}\eta_{ab}\otimes\eta_{cd}
+q^{-2N-1}\eta_{ia}\otimes\eta_{aj}+} \nonumber\\
& &+\delta_{ij}(2Q+Q\s\alpha^{-1}\beta)q^{-2a}\eta_{ab}\otimes
\eta_{ba}+\delta_{ij}(q^{2N+1}Q^{2}-Q\alpha^{-1}\beta)\eta\otimes
\eta).\nonumber
\end{eqnarray}
\end{thm}

\begin{bew}{}
By Lemma \ref{l-invmetrika}, any invariant metric $g$ is of the form
$g(\eta_{ij}\otimes\eta_{kl})=q^{2j}\alpha\delta_{il}\delta_{jk}
+\beta\delta_{ij}\delta_{kl}$ with $\alpha,\beta\in\mathbb{C},
\alpha\not=0,\,\alpha+\s\beta\not=0$. {}From the proof
of Lemma \ref{l-zusa} we know that a connection $\nabla$ on $\Gamma$
is bicovariant if and only if there are complex numbers $\lambda_{1},
\ldots,\lambda_{6}$ such that
\begin{equation} \label{eq-azusansatz}
\nabla(\eta_{ij})=\sum_{n=1}^{6}\lambda_{n}(A_{n}){^{abcd}_{ij}}\eta_{ab}
\otimes\eta_{cd},
\end{equation}
where $\{A_{1},\ldots,A_{6}\}$ generates the vector space
$\mbox{Mor}(v,v\otimes v)$. For our calculi we have $v=u\cont\otimes u$.
By explicit decompositions of the tensor product representations
$u\cont\otimes u$ and $u\cont\otimes u\otimes u\cont\otimes u$ it can be
shown that the following 6 morphisms $A_{1},\ldots,A_{6}$
$(A_{k}=(A_{k}{^{abcd}_{ij}}))$ span the vector space
$\mbox{Mor}(u\cont\otimes u,u\cont\otimes u\otimes u\cont\otimes u)$:
\[
\begin{array}{lll}
A_{1}=q^{-2a-2c}\delta_{ij}\delta^{ab}\delta^{cd} &
A_{2}=q^{-2a}\delta_{ij}\delta^{ad}\delta^{bc} &
A_{3}=\delta{^{a}_{i}}\delta{^{b}_{j}}q^{-2c}\delta^{cd}\\
A_{4}=q^{-2a}\delta_{ab}\delta{^{c}_{i}}\delta{^{d}_{j}} &
A_{5}=\delta{^{a}_{i}}\delta^{bc}\delta{^{d}_{j}} &
A_{6}=q^{2n-2a-2d}\Rda{^{bn}_{am}}\Rda{^{dp}_{cn}}\Rda{^{im}_{jp}}.
\end{array}
\]
For our differential calculus we have
$\dif\eta_{ij}=\eta\wedge\eta_{ij}+\eta_{ij}\wedge\eta$.
Since $\Gamma^{2}=(\Gamma\otimes_{\cal A}\Gamma)/\ker(\sigma-\mbox{id})$,
the torsion of $\nabla$ vanishes if and only if we have
\[
(\sigma-\mbox{id})(\eta\otimes\eta_{ij}+\eta_{ij}\otimes\eta-
\sum_{k=1}^{6}\lambda_{k}A_{k}{^{mnrs}_{ij}}\eta_{mn}\otimes\eta_{rs})=0
\]
in the tensor product $\Gamma\otimes_{\cal A}\Gamma$. Comparing the
coefficients of basis elements the latter is equivalent to the equations
\begin{equation} \label{eq-torkoefaplus}
\lambda_{3}=\lambda_{4}=1+\frac{\lambda_{5}}{q^{-1}Q}-\frac{(Q^{2}+1)
\lambda_{6}}{Q}\quad\mbox{\ for\ }\Gamma_{+,z},
\end{equation}
\begin{equation} \label{eq-torkoefaminus}
\lambda_{3}=\lambda_{4}=1-\frac{\lambda_{5}}{q^{-2N-1}Q}+
\frac{\lambda_{6}}{Q}\quad\mbox{\ for }\Gamma_{-,z}.
\end{equation}
Using Lemma \ref{l-vertrneu} it follows that $\nabla$ satisfies the
compatibility condition with the metric $g$ if and only if the following
equations are fulfilled:\\
\parbox{14cm}{
\[\lambda_{4}+Q\lambda_{6}=0,\qquad
\lambda_{5}+q^{-1}\lambda_{6}=0,\qquad
\alpha\lambda_{2}+(\alpha+\s\beta)\lambda_{3}=0,\]
\[(\alpha+\s\beta)\lambda_{1}+\beta\lambda_{2}+
\beta\lambda_{4}+q^{2N}Q(\beta+qQ\alpha)\lambda_{6}=0\]
} \hfill
\parbox{8mm}{\begin{eqnarray} \label{eq-metkoefaplus} \end{eqnarray}}

for $\Gamma_{+,z}$ and\\
\parbox{142mm}{
\[\lambda_{4}=0,\qquad\lambda_{5}+q^{-2N-1}\lambda_{6}=0,\]
\[(\alpha+\s\beta)\lambda_{1}+\beta\lambda_{2}+
(q^{2N+1}Q^{2}\alpha+q^{2N}Q\beta)\lambda_{6}=0,\]
\[\alpha\lambda_{2}+(\alpha+\s\beta)\lambda_{3}+
(-q^{2N+1}Q\alpha+\beta)\lambda_{5}+
(Q\alpha+q^{-1}\beta)\lambda_{6}=0\]
} \hfill
\parbox{8mm}{\begin{eqnarray} \label{eq-metkoefaminus} \end{eqnarray}}

for $\Gamma_{-,z}$.
Some straightforward computations show that the equations
(\ref{eq-torkoefaplus}) and (\ref{eq-metkoefaplus}) resp.\
(\ref{eq-torkoefaminus}) and (\ref{eq-metkoefaminus}) have unique
solutions $\lambda_{i}$ for $\Gamma_{+,z}$ resp.\ $\Gamma_{-,z}$.
Inserting these solutions into (\ref{eq-azusansatz}) we obtain the
formulas given in the theorem.

It is not difficult to check that the connection defined by the above
formulas fulfills all conditions for a Levi-Civita connection.
\end{bew}

Let us look closer at the Levi-Civita connection on $\Gamma=
\Gamma_{+,z}$.
Using Theorem \ref{t-lcza} for the left and right invariant element
$\eta\in\Gammalinv$ we compute $\nabla(\eta)=\frac{Q^{2}}{2\alpha}
(\alpha+\s\beta)(\s q^{-2a}\eta_{ab}\otimes\eta_{ba}-
\eta\otimes\eta)$. Transforming the basis we obtain
\begin{equation} \label{eq-lczaom0}
\nabla(\omega^{0})=\frac{\mu_{+,z}}{\s}\nabla(\eta)=
\frac{z^{2}q^{2}\mu_{+,z}(\alpha+\s\beta)}{2\alpha}
q^{-2a}\omega^{1}_{ab}\otimes\omega^{1}_{ba}.
\end{equation}
Moreover,
$\nabla(\eta_{ij}-\delta_{ij}\s^{-1}\eta)=\frac{Q}{2}(
q^{2n-2a-2d}\Rda{^{bn}_{am}}\Rdam{^{dp}_{cn}}\Rda{^{im}_{jp}}-
q^{-1}\delta{^{a}_{i}}\delta^{bc}\delta{^{d}_{j}}-
Q\delta^{ab}\delta^{c}_{i}\delta{^{d}_{j}}+\\
Q\s^{-1}q^{-2a-2c}\delta_{ij}\delta^{ab}\delta^{cd})\eta_{ab}\otimes
\eta_{cd}$. Transforming the basis once again we get
\begin{eqnarray} \label{eq-lczaom1}
\lefteqn{ \nabla(\omega^{1}_{ij})=\nu_{+,z}\nabla(\eta_{ij}-
\delta_{ij}\s^{-1}\eta)=}\nonumber\\
& &=\frac{1}{2}(zqq^{2n-2a-2d}\Rda{^{bn}_{am}}\Rdam{^{dp}_{cn}}
\Rda{^{im}_{jp}}\omega^{1}_{ab}\otimes\omega^{1}_{cd}-\nonumber\\
& &-z\omega^{1}_{ia}\otimes\omega^{1}_{aj}-Q^{2}\mu_{+,z}^{-1}
(\omega^{1}_{ij}\otimes\omega^{0}+\omega^{0}\otimes\omega^{1}_{ij})).
\end{eqnarray}

What happens with the Levi-Civita connection in the classical limit?
As explained in Section \ref{sec-jurcokalkule}, we consider the
classical limit in the sense that $z\to1$ and $q\to1$. Retaining
the notation introduced in Sections \ref{sec-jurcokalkule} and
\ref{sec-metrik}, formulas
(\ref{eq-lczaom0}) and (\ref{eq-lczaom1}) show that the limit of the
Levi-Civita connection exists and takes the form
\begin{equation}
\nabla^{\rm cl}(\gwomega^{0})=\lim_{q\to1}\nabla(\omega^{0})=
\frac{(N^{2}-1)c_{0}}{4Nc_{1}}\gwomega^{1}_{ab}\otimes\gwomega^{1}_{ba},
\end{equation}
\begin{equation}
\nabla^{\rm cl}(\gwomega^{1}_{ij})=\frac{1}{2}\left(\gwomega^{1}_{aj}
\otimes\gwomega^{1}_{ia}-\gwomega^{1}_{ia}\otimes\gwomega^{1}_{aj}-
\frac{2N}{N^{2}-1}(\gwomega^{1}_{ij}\otimes\gwomega^{0}+\gwomega^{0}
\otimes\gwomega^{1}_{ij})\right).
\end{equation}

\section{Levi-Civita connections on $O_{q}(N)$ and $Sp_{q}(N)$}
\label{sec-lczbcd}

Now we turn to the differential calculi $\Gamma_{+}$ and $\Gamma_{-}$
on the quantum groups $O_{q}(N)$ and $Sp_{q}(N)$. We begin with

\begin{lemma} \label{l-zusbcd}
We have $\rm{dim\ BC}(\Gamma_{\pm})=14$ for $Sp_{q}(4)$ and
$\rm{dim\ BC}(\Gamma_{\pm})=15$ for $O_{q}(N)$, $N\geq3$ and $Sp_{q}(N)$,
$N>4$.
\end{lemma}

\begin{bew}{}
The proof is similar to the proof of Lemma \ref{l-zusa}. In the present
cases ${\cal A}=O_{q}(N)$ and ${\cal A}=Sp_{q}(N)$ we obtain the
decompositions of the tensor product representations
$u\cont\otimes u=[0]\oplus[2]\oplus[1,1]$,
$u\cont\otimes u\otimes u\cont\otimes u=3[0]\oplus6[2]\oplus k[1,1]\oplus$
other terms, where $k=5$ for $Sp_{q}(4)$ and $k=6$ for $O_{q}(N)$,
$N\geq3$ and $Sp_{q}(N)$, $N>4$.
\end{bew}

Let $\Gamma^{\wedge}=\bigoplus^{\infty}_{n=0}\Gamma^{n}$ be a differential
calculus over ${\cal A}$ such that $\Gamma=\Gamma^{1}=\Gamma_{\pm}$.
In contrast to the $SL_{q}(N)$ case we assume that
$\Gamma^{2}=(\Gamma\otimes_{\cal A}\Gamma)/{\cal K}$, where
${\cal K}=\ker(\sigma-\mbox{id})\oplus\ker(\sigma-q^{N}\mbox{id})\oplus
\ker(\sigma-q^{-N}\mbox{id})$ (see the last remarks in Section
\ref{sec-jurcokalkule}). Such an assumption for the higher order
calculus has already been used in \cite{a-CSchWW1}.

In case of the ``ordinary'' classical differential calculus the
2-forms are the quotient of the tensor product of 1-forms by the
eigenspace of the flip operator with eigenvalue 1. Since the eigenvalues
$q^{N}$ and $q^{-N}$ of the braiding map $\sigma$ tend to 1 when
$q\to 1$, the assumption
$\Gamma^{2}=(\Gamma\otimes_{\cal A}\Gamma)/{\cal K}$ means that the
higher order calculus $\Gamma^{\wedge}$ is some sense nearer to the
corresponding construction in the classical case. Moreover, this
assumption is essential in order to prove the following

\begin{thm} \label{t-lczbcd}
Suppose that $g$ is an invariant metric.
There is exactly one Levi-Civita connection on $\Gamma$ with respect to
$g$.
\end{thm}

\begin{bew}{}
The proof is similar to that of Theorem \ref{t-lcza}.
In Lemma \ref{l-invmetrikbcd} we proved that all invariant metrics have
the form $g(\eta_{ij}\otimes\eta_{kl})=
(\alpha B_{\rm 14}B_{\rm 23}C{^{\rm t}_{\rm 1}}
C{^{\rm t}_{\rm 3}}+\beta B_{\rm 12}B_{\rm 34}\Rda_{\rm 23}
C{^{\rm t}_{\rm 1}}C{^{\rm t}_{\rm 3}}+\gamma B_{\rm 12}B_{\rm 34}
C{^{\rm t}_{\rm 1}}C{^{\rm t}_{\rm 3}})_{ijkl}$ with
$\alpha+\z\beta+\s\gamma\not= 0$, $\alpha\not=q\beta$,
$\alpha\not=-q^{-1}\beta$. By Lemma \ref{l-zusbcd}, the dimension of the
vector space $\mbox{Mor}(u\cont\otimes u,u\cont\otimes u\otimes
u\cont\otimes u)$ is at most 15. A closer investigation of the proof
of Lemma \ref{l-zusbcd} shows that the following 15 morphisms
$A_{k}=(A_{k}{^{mnrs}_{ij}})$ generate the vector space.
(In case $N\geq5$ they form a basis of this space.)
\[
\begin{array}{lll}
A_{1}=B{^{\rm t}_{\rm 1}}B{^{\rm t}_{\rm 3}}C_{\rm 12}C_{\rm 34}
B_{\rm 12}C{^{\rm t}_{\rm 1}},&
A_{2}=B{^{\rm t}_{\rm 1}}B{^{\rm t}_{\rm 3}}\Rdam_{\rm 23}C_{\rm 12}
C_{\rm 34}B_{\rm 12}C{^{\rm t}_{\rm 1}},&
A_{3}=B{^{\rm t}_{\rm 1}}B{^{\rm t}_{\rm 3}}C_{\rm 23}C_{\rm 14}
B_{\rm 12}C{^{\rm t}_{\rm 1}},\\
A_{4}=B{^{\rm t}_{\rm 3}}C_{\rm 34},&
A_{5}=B{^{\rm t}_{\rm 1}}B{^{\rm t}_{\rm 3}}\Rda_{\rm 12}C_{\rm 34}
C{^{\rm t}_{\rm 1}},&
A_{6}=B{^{\rm t}_{\rm 1}}C_{\rm 12},\\
A_{7}=B{^{\rm t}_{\rm 1}}B{^{\rm t}_{\rm 3}}C_{\rm 12}\Rda_{\rm 34}
C{^{\rm t}_{\rm 1}},&
A_{8}=B{^{\rm t}_{\rm 3}}C_{\rm 23},&
A_{9}=B{^{\rm t}_{\rm 1}}B{^{\rm t}_{\rm 3}}C_{\rm 23}\Rda_{\rm 14}
C{^{\rm t}_{\rm 1}},\\
A_{10}=B{^{\rm t}_{\rm 1}}B{^{\rm t}_{\rm 3}}\Rda_{\rm 12}C_{\rm 23}
C{^{\rm t}_{\rm 1}},&
A_{11}=B{^{\rm t}_{\rm 1}}B{^{\rm t}_{\rm 3}}\Rda_{\rm 12}C_{\rm 23}
\Rdam_{\rm 14}C{^{\rm t}_{\rm 1}},&
A_{12}=B{^{\rm t}_{\rm 3}}\Rdam_{\rm 23}C_{\rm 34},\\
A_{13}=B{^{\rm t}_{\rm 1}}B{^{\rm t}_{\rm 3}}\Rdam_{\rm 23}C_{\rm 34}
\Rda_{\rm 12}C{^{\rm t}_{\rm 1}},&
A_{14}=B{^{\rm t}_{\rm 1}}B{^{\rm t}_{\rm 3}}\Rda_{\rm 12}\Rda_{\rm 23}
C_{\rm 34}C{^{\rm t}_{\rm 1}},&
A_{15}=B{^{\rm t}_{\rm 1}}B{^{\rm t}_{\rm 3}}\Rda_{\rm 12}\Rda_{\rm 23}
C_{\rm 34}\Rdam_{\rm 12}C{^{\rm t}_{\rm 1}}.
\end{array}
\]
Therefore, by Lemma \ref{l-Dmorph}, we can make the following ansatz
for our Levi-Civita connection $\nabla$:
\[
\nabla(\eta_{ij})=\sum_{k=1}^{15}\lambda_{k}A_{k}{^{mnrs}_{ij}}\eta_{mn}
\otimes\eta_{rs},\qquad\lambda_{k}\in\mathbb{C}.
\]
The condition for the vanishing torsion takes the form
\[
(\sigma-\mbox{id})(\sigma-q^{N}\mbox{id})(\sigma-q^{-N}\mbox{id})
(\eta\otimes\eta_{ij}+\eta_{ij}\otimes\eta-\sum_{k=1}^{15}\lambda_{k}
A_{k}{^{mnrs}_{ij}}\eta_{mn}\otimes\eta_{rs})=0.
\]
This leads to the equations
\[\lambda_{4}=\lambda_{6}-Q\lambda_{15},\quad
\lambda_{5}=\lambda_{7}-Q\lambda_{14},\quad
\lambda_{8}=\z Q\lambda_{6}+\z\lambda_{15}-\z Q,\quad
\lambda_{9}=\z Q\lambda_{7}+\z\lambda_{14},\]
\[\lambda_{11}-q^{-N}\z\lambda_{10}=Q(q^{-N}\z\lambda_{6}+
q^{-N}(-q^{-N}\z^{2}+1)\lambda_{7}+\lambda_{15}-q^{-N}\z-
\frac{1}{\z^{2}+q^{N}}(\z\lambda_{6}+\lambda_{7}+\z\lambda_{14}+
\z^{2}\lambda_{15}-\z)),\]
\[\lambda_{13}=-Q\lambda_{7}+q^{N}\z^{-1}\lambda_{12}-q^{N}Q\lambda_{15}
+\frac{q^{N}}{\z^{2}+q^{N}}
(\z Q\lambda_{6}+Q\lambda_{7}+\z Q\lambda_{14}+\z^{2}Q\lambda_{15}-\z Q).\]
By Lemma \ref{l-vertrneu}, the connection is compatible with the metric if
and only if the following equations are satisfied:
\[
\begin{array}{l}
\lambda_{6}\beta+\lambda_{7}\alpha=0,\\
\lambda_{6}(\alpha-Q\beta)+\lambda_{7}\beta=0,\\
\lambda_{1}(\alpha+\z\beta+\s\gamma)+\lambda_{2}\z^{-1}\gamma+
\lambda_{3}(Q\beta+\gamma)=0,\\
\lambda_{2}\alpha+\lambda_{3}\beta+\lambda_{5}(\alpha+\z\beta
+\s\gamma)+\lambda_{9}(Q\beta+\gamma)+\lambda_{10}\gamma+
\lambda_{13}\z^{-1}\gamma+\lambda_{14}\z(Q\beta+\gamma)=0,\\
\lambda_{2}\beta+\lambda_{3}(\alpha-Q\beta)+\lambda_{4}(\alpha
+\z\beta+\s\gamma)+\lambda_{8}(Q\beta+\gamma)+\lambda_{11}
\gamma+\lambda_{12}\z^{-1}\gamma+\lambda_{15}\z(Q\beta+\gamma)=0,\\
\lambda_{8}(\alpha-Q\beta)+\lambda_{11}\z\beta+\lambda_{12}\beta+
\lambda_{15}\z(\alpha-Q\beta)=0,\\
\lambda_{8}\beta+\lambda_{10}\z\beta+\lambda_{12}\alpha+
\lambda_{14}\z(\alpha-Q\beta)=0,\\
\lambda_{9}(\alpha-Q\beta)+\lambda_{11}\z\alpha+\lambda_{13}\beta+
\lambda_{15}\z\beta=0,\\
\lambda_{9}\beta+\lambda_{10}\z\alpha+\lambda_{13}\alpha+
\lambda_{14}\z\beta=0.
\end{array}
\]
Set $\tilde{\z}=(q^{N}+1)^{-1}\alpha_{1}^{-1}(\z-\z^{-3}q^{3N})$. Some
computations show that the above system of equations admits a unique
solution
\begin{equation}
\begin{array}{lll}
2\lambda_{4}=-Q^{2}(1+\tilde{\z}\beta),&
2\lambda_{5}=-Q^{2}\tilde{\z}\alpha,&
\lambda_{6}=\lambda_{7}=0,\\
2\lambda_{8}=\z Q(-1+\tilde{\z}\beta),&
2\lambda_{9}=\z Q\tilde{\z}\alpha,&
2\lambda_{10}=Q(1-\tilde{\z}\beta),\\
2\lambda_{11}=-Q\tilde{\z}(\alpha-Q\beta),&
2\lambda_{12}=-\z Q\tilde{\z}(\alpha-Q\beta),&
2\lambda_{13}=-\z Q(1+\tilde{\z}\beta),\\
2\lambda_{14}=Q\tilde{\z}\alpha,&
2\lambda_{15}=Q(1+\tilde{\z}\beta),
\end{array}
\end{equation}
\[\lambda_{3}=\frac{Q^{2}\alpha_{0}\alpha}{2\alpha_{1}\alpha_{2}}+
Q\tilde{\z}\gamma,\quad
\lambda_{2}=-\frac{Q^{2}\alpha_{0}\beta}{2\alpha_{1}\alpha_{2}}+
\frac{Q^{2}}{2}\tilde{\z}\alpha_{0}-Q\z\tilde{\z}(Q\beta+\gamma),\]
\[\lambda_{1}=\frac{Q^{2}}{2}\left(\frac{-Q\alpha\beta-\alpha\gamma+
\z^{-1}\beta\gamma}{\alpha_{1}\alpha_{2}}-
\z^{-1}\tilde{\z}\gamma\right).\]
{}From the preceding considerations it is clear that the
corresponding connection $\nabla$ is indeed a Levi-Civita connection
for $g$.
\end{bew}

In order to examine the classical limit, we have to rewrite the
Levi-Civita connection from Theorem \ref{t-lczbcd} in terms of the
standard basis. Using the projections $P_{0}$, $P_{1}$ and $P_{2}$
defined in Section \ref{sec-metrik} some straightforward
computations yield the formulas
\begin{eqnarray*}
\lefteqn{\nabla(\eta)=\frac{Q\alpha_{0}}{2\alpha_{1}}
(\z q^{-N}+\z^{-1}q^{N})(\z^{2}q^{-N}-\z^{-2}q^{N})
B^{\rm t}{^{a}_{m}}C^{md}P_{1}\eta_{ab}\otimes P_{1}\eta_{bd}+}\\
& &+\frac{Q\alpha_{0}}{2}\left(\frac{Q\s}{\alpha_{2}}+\tilde{\z}
(1-q^{N})(1-q^{N}\z^{-2})\right)B^{\rm t}{^{a}_{m}}C^{md}
P_{2}\eta_{ab}\otimes P_{2}\eta_{bd},
\end{eqnarray*}
\begin{eqnarray*}
\lefteqn{\nabla(P_{1}\eta_{ij})=
\frac{Q(\z q^{-N}+\z^{-1}q^{N})}{2}\left(\frac{\z^{-2}q^{N}-\z^{2}q^{-N}}
{\s}(\eta\otimes P_{1}\eta_{ij}+P_{1}\eta_{ij}\otimes\eta)\right.-}\\
& &-(\z^{2}q^{-N}+1)P_{1}{^{kl}_{ij}}P_{1}\eta_{km}\otimes
P_{1}\eta_{ml}+
(\z^{-2}q^{2N}-\z^{2}q^{-N})P_{1}{^{kl}_{ij}}P_{2}\eta_{km}\otimes
P_{1}\eta_{ml}-\\
& &-\z\tilde{\z}\alpha_{1}(1+q^{-N})P_{1}{^{kl}_{ij}}P_{1}\eta_{km}
\otimes P_{2}\eta_{ml}+
\left.\tilde{\z}\alpha_{1}(\z^{-1}q^{N}-\z)P_{1}{^{kl}_{ij}}
P_{2}\eta_{km}\otimes P_{2}\eta_{ml}\right),
\end{eqnarray*}
{\setlength{\arraycolsep}{0pt}
\begin{eqnarray*}
\nabla(P_{2}\eta_{ij})=&&
\frac{Q(\z q^{-N}\!+\!\z^{-1}q^{N})}{2}\left(\frac{q^{-N}\!-\!q^{N}}{\s}\eta
\otimes P_{2}\eta_{ij}+
\frac{(q^{N}\!-\!1)(1\!-\!\z^{-2}q^{N})\tilde{\z}\alpha_{2}\!-\!Q\s}
{(\z q^{-N}\!+\!\z^{-1}q^{N})\s}P_{2}\eta_{ij}\otimes\eta\right.-\\
& &-(q^{N}+\z^{2}q^{-2N})P_{2}{^{kl}_{ij}}P_{1}\eta_{km}\otimes
P_{2}\eta_{ml}+(1-q^{N})P_{2}{^{kl}_{ij}}P_{2}\eta_{km}\otimes
P_{2}\eta_{ml}+\\
& &+\z\tilde{\z}\alpha_{2}(1+q^{-N})P_{2}{^{kl}_{ij}}P_{1}\eta_{km}
\otimes P_{1}\eta_{ml}+
\left.\tilde{\z}\alpha_{2}(\z-\z^{-1}q^{N})P_{2}{^{kl}_{ij}}
P_{2}\eta_{km}\otimes P_{1}\eta_{ml}\right).
\end{eqnarray*}
}

{}From these formulas it follows that the classical limits of the
Levi-Civita connections $\nabla$ for both calculi $\Gamma_{+}$ and
$\Gamma_{-}$ exist. For the subspaces $\Upsilon^{1}$ of 1-forms of
$\Gamma_{+}$ which corresponds to the classical differential calculus
(see \cite{a-HSS1}) we obtain
\begin{eqnarray*}
\lefteqn{\nabla^{\rm cl}(\gwomega^{1}_{ij})=\lim_{q\to1}\nabla
(\omega^{1}_{ij})=-\epsilon\frac{N-2\epsilon}{N-\epsilon}
(\gwomega^{0}\otimes\gwomega^{1}_{ij}+\gwomega^{1}_{ij}\otimes
\gwomega^{0})-\gwP_{1}{^{kl}_{ij}}\gwomega^{1}_{km}\otimes
\gwomega^{1}_{ml}-}\\
&&-\frac{N-4\epsilon}{N-2\epsilon}
\gwP_{1}{^{kl}_{ij}}(\gwomega^{1}_{km}\otimes\gwomega^{2}_{ml}+
\gwomega^{2}_{km}\otimes\gwomega^{1}_{ml}+\gwomega^{2}_{km}\otimes
\gwomega^{2}_{ml}).
\end{eqnarray*}

\begin{appendix}

\section{The Rosso form of ${\cal U}_{q}(sl(N))$}
\label{sec-Rossoform}

In Section \ref{sec-metrik} we defined invariant metrics and we have
seen that such metrics are not uniquely determined. On the other hand,
Rosso showed in \cite{a-Rosso2} that there is a unique ad-invariant
bilinear form for the quantum universal enveloping algebra
${\cal U}_{q}(\mathfrak{g})$ for a simple Lie algebra $\mathfrak{g}$.
In this appendix we define such a form adapted to the
preceding considerations and we compute the corresponding ad-invariant
metric on the quantum Lie algebra ${\cal X}$ of the FODC
$\Gamma_{+,z}$.

Let $q$ be a complex number, $q\not=0$, $q^{k}\not=1$ for all
$k\in\mathbb{N}$ and let $(a_{ij})$ be the Cartan-matrix for $sl(N)$.
Let $({\cal U}_{q}(sl(N)),\Delta,\kappa,\varepsilon)$ be the Hopf algebra
over $\mathbb{C}$ generated by the set of elements
$\{E_{i},F_{i},K_{i},K_{i}^{-1},\tilde{K}_{i},\tilde{K}_{N},
\tilde{K}_{i}^{-1},\tilde{K}_{N}^{-1}\,|\,i=1,\ldots,N-1\}$
with relations\\
\parbox{14cm}{
\[K_{i}K_{j}=K_{j}K_{i},\quad K_{i}K_{i}^{-1}=K_{i}^{-1}K_{i}=1,
\quad K_{i}\tilde{K}_{n}=\tilde{K}_{n}K_{i},\]
\[\tilde{K}_{m}\tilde{K}_{n}=\tilde{K}_{n}\tilde{K}_{m},\quad
\tilde{K}_{n}\tilde{K}_{n}^{-1}=\tilde{K}_{n}^{-1}\tilde{K}_{n}=1,
\quad \tilde{K}_{i}^{-1}\tilde{K}_{i+1}=K_{i}^{2},\]
\[K_{i}E_{j}=q^{a_{ij}}E_{j}K_{i},\qquad
K_{i}F_{j}=q^{-a_{ij}}F_{j}K_{i}\]
\[\tilde{K}_{n}E_{j}=
q^{2(\delta_{n,j+1}-\delta_{nj})}E_{j}\tilde{K}_{n}\qquad
\tilde{K}_{n}F_{j}=
q^{2(\delta_{nj}-\delta_{n,j+1})}F_{j}\tilde{K}_{n},\]
\[E_{i}F_{j}-q^{-a_{ij}}F_{j}E_{i}=
\delta_{ij}\frac{K_{i}^{2}-1}{q^{2}-1},\]
\[\sum_{k=0}^{1-a_{ij}}(-1)^{k}b_{-a_{ij},k}
E_{i}^{k}E_{j}E_{i}^{1-a_{ij}-k}=0\quad(i\not= j),\]
\[\sum_{k=0}^{1-a_{ij}}(-1)^{k}b_{-a_{ij},k}
F_{i}^{k}F_{j}F_{i}^{1-a_{ij}-k}=0\quad(i\not= j)\]
} \hfill
\parbox{8mm}{\begin{equation} \label{eq-uqslnalg} \end{equation}}\\
and coproduct $\Delta$, antipode $\kappa$, counit $\varepsilon$
defined by\\
\parbox{14cm}{
\[\Delta(E_{i})=E_{i}\otimes K_{i}+1\otimes E_{i},\qquad
\Delta(F_{i})=F_{i}\otimes K_{i}+1\otimes F_{i},\]
\[\Delta(K_{i})=K_{i}\otimes K_{i},\qquad
\Delta(\tilde{K}_{n})=\tilde{K}_{n}\otimes\tilde{K}_{n},\]
\[\kappa(E_{i})=-q^{2}K_{i}^{-1}E_{i},\quad
\kappa(F_{i})=-F_{i}K_{i}^{-1},\quad
\kappa(K_{i})=K_{i}^{-1},\quad
\kappa(\tilde{K}_{n})=\tilde{K}_{n}^{-1},\]
\[\varepsilon(E_{i})=\varepsilon(F_{i})=0,\quad
\varepsilon(K_{i})=\varepsilon(\tilde{K}_{n})=1\]
} \hfill
\parbox{8mm}{\begin{equation} \label{eq-uqslnhopf} \end{equation}}\\
for all $i,j=1,\ldots,N-1$ and all $m,n=1,\ldots,N$. The constants
in (\ref{eq-uqslnalg}) are $b_{0,0}=b_{0,1}=b_{1,0}=b_{1,2}=1$,
$b_{1,1}=q+q^{-1}$. A realization of the Hopf algebra
${\cal U}_{q}(sl(N))$ in terms of the $L$-functionals is obtained by
setting
\begin{equation} \label{eq-uqgener}
E_{i}=Q^{-1}l^{-}{^{i}_{i}}l^{+}{^{i}_{i+1}},\quad
F_{i}=-q^{-1}Q^{-1}l^{-}{^{i+1}_{i}}l^{+}{^{i+1}_{i+1}},\quad
K_{i}=l^{-}{^{i}_{i}}l^{+}{^{i+1}_{i+1}},\quad
\tilde{K}_{i}=(l^{+}{^{i}_{i}})^{2}.
\end{equation}

This Hopf algebra is $\mathbb{Z}^{N-1}$-graduated with grading
$\partial$ given by $\partial E_{i}:=\alpha_{i}$,
$\partial F_{i}=-\alpha_{i}$, $\partial K_{i}=0$,
$\partial\tilde{K}_{i}=0$. Let $<\cdot,\cdot>$
denote the symmetric bilinear form such that
$<\alpha_{i},\alpha_{j}>=a_{ij}$.

It is well-known that all elements of ${\cal U}_{q}(sl(N))$ can be written
as finite linear combinations of terms of the form $FKE$, where $F$
and $E$ are finite products of elements $F_{i}$ and $E_{i}$,
respectively, and $K$ is a product of the generators
$K_{i}$, $K_{i}^{-1}$, $\tilde{K}_{i}$ and $\tilde{K}_{i}^{-1}$.

Let us recall that a map $(\cdot,\cdot):{\cal B}\times{\cal B}\to
\mathbb{C}$ for a Hopf algebra ${\cal B}$ is called ad-invariant
if
\begin{equation}\label{eq-adinvmap}
(\adR\xi_{(1)}(\zeta_{1}),\adR\xi_{(2)}(\zeta_{2}))=
\varepsilon(\xi)(\zeta_{1},\zeta_{2})\quad
\text{for $\xi,\zeta_{1},\zeta_{2}\in{\cal B}$}.
\end{equation}
Then one can prove the following proposition which is essentially
Rosso's result adapted to the present setting.
Note that in \cite{a-Rosso2} the algebra ${\cal U}_q(sl(N))$ is different
from ours and the left adjoint action is used.

\begin{satz} \label{s-rossomap}
There is a unique ad-invariant bilinear map
$(\cdot,\cdot):{\cal U}_{q}(sl(N))\times{\cal U}_{q}(sl(N))
\to\mathbb{C}$
such that $(FKE,F'K'E')=(F,E')(K,K')(E,F')$, $(K,K')=(K',K)$,
$(KK',K'')=$\\
$(K,K'')(K',K'')$, $(K_{i},K_{j})=q^{-a_{ij}/2}$,
$(\tilde{K}_{i},K_{j})=q^{\delta_{ij}-\delta_{i,j+1}}$
and $(\tilde{K}_{i},\tilde{K}_{j})=q^{-2\delta_{ij}}$.
\end{satz}

\begin{bew}{}
The proof is similar to that of Theorem 6 in \cite{a-Rosso2}. We omit
the details. Howeover, we want to stress that we deal with another
adjoint action and with different commutation relations of the
generators of ${\cal U}_{q}(sl(N))$.
\end{bew}

In what follows we use the abbreviations
$E_{i,i+1}:=E_{i}$, $F_{i+1,i}:=F_{i}$,
$E_{i,j+1}:=E_{i+1,j+1}E_{i}-q^{-1}E_{i}E_{i+1,j+1}$,
$F_{j+1,i}:=F_{i}F_{j+1,i+1}-q^{-1}F_{j+1,i+1}F_{i}$ for $i<j$.
Using (\ref{eq-adinvmap}) and Proposition \ref{s-rossomap}
the bilinear form $(\cdot,\cdot)$ for these elements of
${\cal U}_{q}(sl(N))$ can be computed. The result is given by the
formulas
\begin{eqnarray*}
(E_{i},F_{j})=-q^{-1}Q^{-1}\delta_{ij},&&
(E_{ij},F_{kl})=-q^{2i-2j+1}Q^{-1}\delta_{il}\delta_{jk},\\
(F_{i},E_{j})=-qQ^{-1}\delta_{ij},&&
(F_{ij},E_{kl})=-qQ^{-1}\delta_{il}\delta_{jk}.
\end{eqnarray*}
{}From equations (\ref{eq-chifunkt}) and (\ref{eq-uqgener}) the
generators $\chi_{ij}$ of the quantum Lie algebra ${\cal X}$
of the FODC $\Gamma_{+,z}$ can be expressed in terms of the
elements $E_{ij}$, $F_{ji}$, $K_{i}$ and $\tilde{K}_{i}$:
\begin{eqnarray*}
\chi_{ij}&=&q^{-1}QF_{ji}\tilde{K}_{i}+q^{-1}Q^{2}\sum_{r<i}
F_{jr}\tilde{K}_{r}E_{ri},\quad (i<j)\\
\chi_{ij}&=&q^{-1}Q\tilde{K}_{j}E_{ji}+q^{-1}Q^{2}\sum_{r<j}
F_{jr}\tilde{K}_{r}E_{ri},\quad (i>j)\\
\chi_{ii}&=&q^{-2}\tilde{K}_{i}-q^{-1}Q\sum_{m=1}^{i-1}q^{2m-2i}
\tilde{K}_{m}+q^{-2}Q^{2}\sum_{n<i}F_{in}\tilde{K}_{n}E_{ni}-\\
&&{}-q^{-1}Q^{3}\sum_{m=1}^{i-1}\sum_{n=1}^{m-1}q^{2m-2i}F_{mn}
\tilde{K}_{n}E_{nm}-q^{-2i}\varepsilon.
\end{eqnarray*}

Combining both sets of the preceding formulas it follows that the
bilinear form in Proposition \ref{s-rossomap} takes the following
form on the generators $\chi_{ij}$ of the quantum Lie algebra:
\begin{equation}
(\chi_{ij},\chi_{kl})=-q^{-2i-1}Q\delta_{il}\delta_{jk}+
q^{-2i-2j}Q^{2}\delta_{ij}\delta_{kl}.
\end{equation}
This is precisely the dual metric $g^{*}$ described by formula
(\ref{eq-dualmeta}) for the parameter values $\alpha=-qQ^{-1}$,
$\beta=-q^{2N+2}$ (cf.\ (\ref{eq-dualmeta})).

\section{Another way to Levi-Civita connections}
\label{sec-LCZ2}

The results in Theorem \ref{t-lcza} and \ref{t-lczbcd} indicate that the
concepts defined above are useful. Nevertheless, there are various
other possibilities to define a metric and the compatibility of a metric
with a connection. Here we pick out one of these possibilities which is
similar to the classical case.

Let $\Gamma$ denote one of the first order differential calculi
$\Gamma_{\pm,z}$ of $SL_{q}(N)$ defined in Section
\ref{sec-jurcokalkule}. Suppose that $q$ is real. Then $\Gamma$ is a
$*$-calculus for the Hopf-$*$-algebra ${\cal A}=SU_{q}(N)$. The
involutions of ${\cal A}$ and of $\Gamma$ are given by
$(u^{i}_{j})^{*}=\kappa(u^{j}_{i})$ and ${\eta_{ij}}^{*}=-\eta_{ji}$,
respectively.

We call a map $g:\Gamma\times\Gamma\rightarrow{\cal A}$ a {\bf metric} if
\[g(a\xi,b\zeta)=ag(\xi,\zeta)b^{*}\qquad\mbox{for all }a,b\in{\cal A},
\quad\xi,\zeta\in\Gamma.\]
We say a (left) connection $\nabla$ on $\Gamma$ is compatible with the
metric $g$ if
\[g(\nabla(\xi),\zeta)+g(\xi,\nabla(\zeta))=\dif g(\xi,\zeta)\qquad
\mbox{for all }\xi,\zeta\in\Gamma.\]
This is well defined because of\\
$g(\nabla(a\xi),b\zeta)+g(a\xi,\nabla(b\zeta))-\dif g(a\xi,b\zeta)=
\dif(a)g(\xi,\zeta)b^{*}+ag(\nabla(\xi),\zeta)b^{*}+
ag(\xi,\zeta)\dif(b^{*})+\\
ag(\xi,\nabla(\zeta))b^{*}-
\dif(a)g(\xi,\zeta)b^{*}-a\dif g(\xi,\zeta)b^{*}-
ag(\xi,\zeta)\dif(b^{*})=a\left(g(\nabla(\xi),\zeta)+
g(\xi,\nabla(\zeta))-\dif g(\xi,\zeta)\right)b^{*}.$

A bicovariant left connection $\nabla$ is called Levi-Civita connection
in respect to the metric $g$, if $\nabla$ has vanishing torsion and if it
is compatible with the metric $g$. As above, such a Levi-Civita
connection is uniquely determined by its values on the generators
$\eta_{ij}\in\Gammalinv$.

Let $g'$ be a metric from Lemma \ref{l-invmetrika} with
$\alpha,\beta\in\mathbb{R}$ and let $g$ defined
by $g(\xi,\zeta)=g'(\xi\otimes\zeta^{*})$ with $\xi\in\Gamma$,
$\zeta\in\Gammalinv$. Then we have
$g(\eta_{ij},\eta_{kl})=q^{2j}\alpha\delta_{ik}\delta_{jl}+
\beta\delta_{ij}\delta_{kl}$, $\alpha\not= 0$,
$\alpha+\s\beta\not= 0$.
Let us assume that $\nabla(\eta_{ij})=\sum_{k=1}^{6}\lambda_{k}
A_{k}{^{mnrs}_{ij}}\eta_{mn}\otimes\eta_{rs}$ as in the proof of Theorem
\ref{t-lcza}. In order to simplify the calculation, let $\lambda_{k}$ be
real.

The torsion of $\nabla$ has to vanish and so we get (as in Theorem
\ref{t-lcza})
\[Q\lambda_{3}=Q\lambda_{4}=q\lambda_{5}-(Q^{2}+1)\lambda_{6}+Q.\]
The compatibility with the metric $g$ gives only one equation
\[\lambda_{2}\alpha-\lambda_{3}(\alpha+\s\beta)-
\lambda_{5}\beta-q^{-1}\lambda_{6}\beta=0.\]
Hence there is a three parameter family of Levi-Civita connections.
Moreover, we have
\[\nabla(\eta)=(\lambda_{2}\s+\lambda_{5}+q^{-1}\lambda_{6})
q^{-2m}\eta_{mn}\otimes\eta_{nm}+
(\lambda_{1}\s+\lambda_{3}+\lambda_{4}+q^{2N}Q\lambda_{6})\eta\otimes\eta
\]
for the bi-invariant element $\eta\in\Gamma$. That is, even if we require
in addition that $\nabla(\eta)=\lambda\eta\otimes\eta$ for some
$\lambda\in\mathbb{C}$ we do not get a {\bf unique} Levi-Civita connection.
\end{appendix}

\nocite{a-GMMM1} \nocite{a-Mourad1} \nocite{b-Kassel1} \nocite{a-Drinfeld1}
\nocite{a-Jimbo1} \nocite{a-BDMS1} \nocite{a-DVM1} \nocite{a-CQ1}
\nocite{a-BrzMaj1} \nocite{a-Majid1}

%\bibliography{quantum}
%\bibliographystyle{mybibeng}

\end{document}